\documentclass[11pt]{article}
\pdfoutput=1
\usepackage{amsfonts,amssymb,amsthm,amsmath,mathrsfs,mathtools}
\usepackage{url}
\expandafter\def\expandafter\UrlBreaks\expandafter{\UrlBreaks
  \do\a\do\b\do\c\do\d\do\e\do\f\do\g\do\h\do\i\do\j%
  \do\k\do\l\do\m\do\n\do\o\do\p\do\q\do\r\do\s\do\t%
  \do\u\do\v\do\w\do\x\do\y\do\z\do\A\do\B\do\C\do\D%
  \do\E\do\F\do\G\do\H\do\I\do\J\do\K\do\L\do\M\do\N%
  \do\O\do\P\do\Q\do\R\do\S\do\T\do\U\do\V\do\W\do\X%
  \do\Y\do\Z}

\newtheorem{theorem}{\bf Theorem}

\newtheorem{definition}{\bf Definition}

\title{Simpson's Paradox: A Singularity of Statistical and Inductive Inference}

\author{
Palash Sarkar\\
Applied Statistics Unit\\
Indian Statistical Institute\\
203, B.T. Road, Kolkata\\
India 700108\\
email: palash@isical.ac.in
\and
Prasanta S. Bandyopadhyay \\
Department of History \& Philosophy \\
Montana State University \\
Bozeman, MT \\
USA \\
email: psb@montana.edu
}

\newcommand{\sym}[1]{{\sf #1}}

\begin{document}

\maketitle



%

\begin{abstract}
The occurrence of Simpson's paradox (SP) in $2\times 2$ contingency tables has been well studied. The present work comprehensively 
	revisits this problem using a combination of philosophical reflections, causal considerations, and probability theory. The first contribution is to provide a 
	schematic analysis of SP in 
	$2\times 2$ contingency tables and present new results, detailed proofs of previous results and a unifying view of the important examples of SP that have been 
	reported in the literature. The second contribution of the paper suggests a new perspective on the surprise element of SP, raises some critical questions regarding 
	the influential causal analyses of SP and provides a broad perspective on logic, probability, and statistics with SP at its center. The upshot of this research 
	is that we need both causal concepts and statistical tools coupled with philosophical analyses to sort out issues regarding SP. 
\end{abstract}

\section{Introduction \label{sec-intro} }
Simpson's paradox (SP) involves the reversal of the direction of a comparison or the cessation of the association when data from several groups are combined. 
Like most paradoxes, SP challenges our intuitions in a fundamental way. We discuss how human intuitions often go astray in the case of SP when intuitions about 
arithmetical sum and ratio diverge radically. SP is critically important in the current multidisciplinary approaches to causal inference. According to Judea Pearl, one 
of the pioneers among causal theorists, causal intuition is the only way to untangle the paradox. Philosophers/statisticians Peter Spirtes, Clark Glymour, and Richard Scheines 
(2000, hereafter the Carnegie Mellon University or CMU group) take one aim of their work to be ``[clarifying] diverse topics,'' including SP. Focusing on 
the paradox will shed light on what the causalists and statisticians of various stripes (Pearson et al. 1899, Yule, 1903, Blyth, 1972, Armistead, 2014, Simpson, 1951, and 
Mittal, 1991) want to say about it, offering a vantage point regarding this vast area of current research. 

Statistical, and more generally inductive, inference is arguably one of the most important methods
of knowledge acquisition. SP with its uniquely distinctive characteristic remains a sigularity in the
theory of statistical and inductive inference. Our work continues the long line of research in elucidating
and understanding SP from both mathematical and philosophical points of view.

Keeping in mind that ``definition has a distinctively philosophical use (Quine, 1975)'', we provide a definition of SP and then provide for the first time an exhaustive 
analysis of necessary and sufficient conditions for SP in a $2\times 2$ contingency table going beyond what is discussed in the Stanford Encyclopedia of Philosophy 
(Sprenger and Weinberger, 2021).  Unless we know how SP could result, it would be difficult to appreciate its different dimensions and the debate over its true nature. Recently, 
there has been a call for pluralism in philosophy in which various perspectives such as in logic, epistemology, and the rest, request a respectful critical appraisal. Formal 
epistemology need not be immune from it. To this end, we propose a non-causal account of the paradox which represents a minority view. 
We will evaluate David Lewis's influential dictum as a minimal constraint any account of causality must satisfy. He writes, ``[w]hen commonsense delivers a firm and 
uncontroversial answer about a not-too-far-fetched case, theory had better agree. If an analysis of causation does not deliver the common-sense answer that is bad trouble.'' 
We will argue that because of the complex nature of SP, we fail to confirm or disconfirm this dictum decisively in appraising any causal theory.
The paradox comes as a surprise because of its reversal or 
cessation of an association in the combined population, which one observes in the sub-populations. Causalists explain this element of paradoxicality in SP solely in causal 
terms and intuitions. We disagree. Our diagnosis of its paradoxicality is provided in terms of the violation of a non-causal principle avowedly rejected by causal theorists 
such as Pearl as the sole reason for its paradoxicality. 

We tend to agree with the view that philosophical discourse is about ``arguments'' rather than about ``proofs'' (Sober, 2021.) Providing proofs, however,
and improving on the existing ones will certainly strengthen a specific philosophical argument, thus lending support for the stance one advocates. 
The literature contains a number of results providing necessary conditions for SP and sufficient conditions for SP not to hold. 
See Sprenger and Weinberger (2021) for a summary of such results.
We comprehensively revisit this line of research and obtain new results as well as new and improved proofs of important prior results. 
Statistical inference is an inductive procedure. With progressive accumulation of data, it is possible that SP can arise and vanish over time. We show that this can indeed 
happen where the appearance and disappearance of SP can toggle ad infinitum. 

A key question when faced with an instance of SP, is whether to act based upon the sub-population data or to act based upon the aggregated data. It has 
been suggested by Pearl (2014) that depending on the situation, the correct answer lies in either the aggregated, or the disaggregated data, or 
neither\footnote{An explanation of Simpson's paradox from the view point of confirmation theory has been put forward by Fitelson (2017). While the explanation is interesting, the work
does not consider the important problem of ``what to do?'' when faced with an instance of SP. In fact, it would be of interest to perform a detailed analysis
of the examples described in Section~\ref{sec-ex} based on Fitelson's approach.}.
We make two points regarding this suggestion. First, that that there can be situations where both the aggregated and the disaggregated data can 
provide meaningful answers, i.e., for the same data and the same context, depending upon extraneous information or upon the question that 
one wishes to be answered, the correct choice is either the aggregated or the disaggregated table. We provide concrete examples in Section~\ref{sec-resolution}. 
So, a priori, one cannot rule out one of these cases.
Second, that among the examples of SP that actually appear in the literature, there is none where the correct answer does not lie in either the 
aggregated or the disaggregated data, raising the question of whether the possibility of not being able to obtain the correct answer from either the aggregated or the
disaggregated data is only a mathematical artifice. 

In the final section, we provide a broad discussion of the interplay of logic, probability and philosophy revolving around causality. The discussion is anchored in 
philosophy of statistics with SP at its core. We put forward the view that an appropriate amalgamation of statistical and causal techniques can lead to improved methods and 
appreciation of scientific inference.

The paper contains a number of theorems. The proofs are provided in the appendix.

\section{Definition of Simpson's Paradox \label{sec-def-SP} }
In this section, we take a detailed look at what constitutes a definition of Simpson's Paradox in terms of $2\times 2$ contingency tables. There are 
three such distinct definitions of SP in the literature, namely a definition by Blyth (1972), a definition by Mittal (1991) and a definition by
Bandyopadhyay et al. (2011). The three definitions do not agree on all cases. In view of this, it is worthwhile to consider the possible cases that can
arise and identify which of these constitute examples of SP. The goal of the present section is to perform this task in a comprehensive manner.

Let $T$ be a $2\times 2$ contingency table of the form given below.
\begin{eqnarray} \label{eqn-cont-table}
T & = & 
\begin{array}{|c|c|}
\hline
a & b \\ \hline
c & d \\ \hline
\end{array}\ .
\end{eqnarray}
We will assume that $a$, $b$, $c$ and $d$ are positive real numbers.

Given a contingency table $T$ as in~\eqref{eqn-cont-table}, one may associate true/false valued random variables
$X$, $Y$ and $\overline{X}$ and $\overline{Y}$ as follows: $X$ and $\overline{X}$ correspond to the first and second row respectively of $T$;
$Y$ and $\overline{Y}$ correspond to the first and second column respectively of $T$. 
Note that $X$ (resp. $Y$) is true if and only if $\overline{X}$ (resp. $\overline{Y}$) is false.
So, $\Pr[Y|X] = a/(a+b)$ and $\Pr[Y|\overline{X}]=c/(c+d)$ and similarly for the other conditional probabilities. 
Here the notation $\Pr[Y|X]$ is implicitly used to denote the probability of the event $\Pr[Y=\mbox{true}|X=\mbox{true}]$ and similarly for other probabilities.

\begin{definition}\label{defn-cont-tab}
Let $T_1$ and $T_2$ be two $2\times 2$ contingency tables as shown below.
\begin{center}
\begin{tabular}{cc}
\begin{minipage}{100pt}
\begin{eqnarray*}
T_1 & = & 
\begin{array}{|c|c|}
\hline
a_1 & b_1 \\ \hline
c_1 & d_1 \\ \hline
\end{array}\ ,
\end{eqnarray*}
\end{minipage}
\begin{minipage}{100pt}
\begin{eqnarray*}
T_2 & = & 
\begin{array}{|c|c|}
\hline
a_2 & b_2 \\ \hline
c_2 & d_2 \\ \hline
\end{array}\ .
\end{eqnarray*}
\end{minipage}
\end{tabular}
\end{center}

The collapsed (or aggregated) table $T_1+T_2$ is the following.
\begin{eqnarray*}
T_1+T_2 & = & 
\begin{array}{|c|c|}
\hline
a_1+a_2 & b_1+b_2 \\ \hline
c_1+c_2 & d_1+d_2 \\ \hline
\end{array}\ .
\end{eqnarray*}
\end{definition}

\paragraph{Notation:} We fix the following notation for later use.
\begin{eqnarray*}
\alpha_1 = a_1+b_1; \quad \alpha_2=a_2+b_2; \quad \gamma_1=c_1+d_1; \quad \gamma_2=c_2+d_2; \\
A_1 = \frac{a_1}{\alpha_1}; \quad A_2 = \frac{a_2}{\alpha_2}; \quad C_1 = \frac{c_1}{\gamma_1}; \quad C_2 = \frac{c_2}{\gamma_2}; \\
	\mu = \frac{a_1+a_2}{\alpha_1+\alpha_2}; \quad \nu = \frac{c_1+c_2}{\gamma_1+\gamma_2}.
\end{eqnarray*}

We connect the above quantities to conditional probabilities. Let $M$ and $\overline{M}$ be true/false valued random variables associated with tables $T_1$ and
$T_2$ respectively; $X$ and $\overline{X}$ (resp. $Y$ and $\overline{Y}$) correspond to the first and second rows (resp. columns) of $T_1,T_2$ and $T_1+T_2$.
Then $\Pr[Y|X,M] = A_1$, $\Pr[Y|\overline{X},M]=C_1$, $\Pr[Y|X,\overline{M}] = A_2$, $\Pr[Y|\overline{X},\overline{M}]=C_2$, $\Pr[Y|X]=\mu$
and $\Pr[Y|\overline{X}]=\nu$.

\begin{figure}
	\caption{Classifications of the 27 relationships that can occur for the tables $T_1$, $T_2$ and $T_1+T_2$. \label{fig-reln} }
\begin{center}
\begin{tabular}{|l|c|c|c|c|}
\hline
	Case & cond on $T_1$ & cond on $T_2$ & cond on $T_1+T_2$ & remark \\ \hline\hline
	1 & $A_1>C_1$ & $A_2>C_2$ & $\mu>\nu$ & Aligned \\ \hline
	2 & $A_1>C_1$ & $A_2>C_2$ & $\mu=\nu$ & Weak paradox \\ \hline
	3 & $A_1>C_1$ & $A_2>C_2$ & $\mu<\nu$ & Paradox \\ \hline
	4 & $A_1>C_1$ & $A_2=C_2$ & $\mu>\nu$ & Class-2 \\ \hline
	5 & $A_1>C_1$ & $A_2=C_2$ & $\mu=\nu$ & Class-3 \\ \hline
	6 & $A_1>C_1$ & $A_2=C_2$ & $\mu<\nu$ & Class-4 \\ \hline
	7 & $A_1>C_1$ & $A_2<C_2$ & $\mu>\nu$ & Class-1 \\ \hline
	8 & $A_1>C_1$ & $A_2<C_2$ & $\mu=\nu$ & Class-5 \\ \hline
	9 & $A_1>C_1$ & $A_2<C_2$ & $\mu<\nu$ & Class-1 \\ \hline
	10 & $A_1=C_1$ & $A_2>C_2$ & $\mu>\nu$ & Class-2 \\ \hline
	11 & $A_1=C_1$ & $A_2>C_2$ & $\mu=\nu$ & Class-3 \\ \hline
	12 & $A_1=C_1$ & $A_2>C_2$ & $\mu<\nu$ & Class-4 \\ \hline
	13 & $A_1=C_1$ & $A_2=C_2$ & $\mu>\nu$ & Class-6 \\ \hline
	14 & $A_1=C_1$ & $A_2=C_2$ & $\mu=\nu$ & Class-0 \\ \hline
	15 & $A_1=C_1$ & $A_2=C_2$ & $\mu<\nu$ & Class-6 \\ \hline
	16 & $A_1=C_1$ & $A_2<C_2$ & $\mu>\nu$ & Class-4 \\ \hline
	17 & $A_1=C_1$ & $A_2<C_2$ & $\mu=\nu$ & Class-3 \\ \hline
	18 & $A_1=C_1$ & $A_2<C_2$ & $\mu<\nu$ & Class-2 \\ \hline
	19 & $A_1<C_1$ & $A_2>C_2$ & $\mu>\nu$ & Class-1 \\ \hline
	20 & $A_1<C_1$ & $A_2>C_2$ & $\mu=\nu$ & Class-5 \\ \hline
	21 & $A_1<C_1$ & $A_2>C_2$ & $\mu<\nu$ & Class-1 \\ \hline
	22 & $A_1<C_1$ & $A_2=C_2$ & $\mu>\nu$ & Class-4 \\ \hline
	23 & $A_1<C_1$ & $A_2=C_2$ & $\mu=\nu$ & Class-3 \\ \hline
	24 & $A_1<C_1$ & $A_2=C_2$ & $\mu<\nu$ & Class-2 \\ \hline
	25 & $A_1<C_1$ & $A_2<C_2$ & $\mu>\nu$ & Paradox \\ \hline
	26 & $A_1<C_1$ & $A_2<C_2$ & $\mu=\nu$ & Weak paradox \\ \hline
	27 & $A_1<C_1$ & $A_2<C_2$ & $\mu<\nu$ & Aligned \\ \hline
\end{tabular}
\end{center}
\end{figure}
From the tables $T_1$, $T_2$ and $T_1+T_2$, there are 27 relationships that arise. These are summarised in Figure~\ref{fig-reln}. 
The last column of the table provides a classification of these relationships. The basis for this classification is as follows.
\begin{description}
	\item{\em Class-0:} Equality holds in all the three tables $T_1$, $T_2$ and $T_1+T_2$.
	\item{\em Aligned:} The inequalities in the Tables $T_1$, $T_2$ and $T_1+T_2$ are strict and are all in the same direction.
	\item{\em Paradox:} The inequalities in the Tables $T_1$, $T_2$ and $T_1+T_2$ are strict; the inequalities in the Tables $T_1$ and $T_2$ are in the same direction, while
		the inequality in the Table $T_1+T_2$ is in the other direction.
	\item{\em Weak paradox:} The inequalities in the Tables $T_1$ and $T_2$ are strict and are in the same direction, while there is equality in Table $T_1+T_2$. 
	\item{\em Class-1:} The inequalities in the Tables $T_1$, $T_2$ and $T_1+T_2$ are strict. The inequalities in $T_1$ and $T_2$ are in opposite directions.
	\item{\em Class-2:} One of the inequalities in Tables $T_1$ and $T_2$ is strict, while the other is a equality. The inequality in Table $T_1+T_2$ is strict
		and is in the same direction as the strict inequality in $T_1$ or $T_2$. 
	\item{\em Class-3:} One of the inequalities in Tables $T_1$ and $T_2$ is strict, while the other is a equality. Equality holds in Table $T_1+T_2$.
	\item{\em Class-4:} One of the inequalities in Tables $T_1$ and $T_2$ is strict, while the other is a equality. The inequality in Table $T_1+T_2$ is strict
		and is in the opposite direction to the strict inequality in $T_1$ or $T_2$. 
	\item{\em Class-5:} The inequalities in the Tables $T_1$ and $T_2$ are strict and in the opposite directions. Equality holds in Table $T_1+T_2$. 
	\item{\em Class-6:} Equality holds in the Tables $T_1$ and $T_2$, while strict inequality holds in Table $T_1+T_2$. 
\end{description}
A relevant question is whether each of the cases in Figure~\ref{fig-reln} is possible, i.e., for each case, is it possible to exhibit a pair of $2\times 2$ contingency tables
which satisfies the conditions of that particular case? We show that the answer to this question is positive. 
\begin{theorem}\label{thm-possible}
All the cases considered in Figure~\ref{fig-reln} are logically possible. 
\end{theorem}	


Equality represents a no-preference condition, while a strict inequality represents a clear preference of one row over the other. Let us try to understand 
the entries in Figure~\ref{fig-reln} in terms of the interpretation of equality as no-preference and strict inequality as a clear preference. 
The understanding of Class-0 and the `Aligned' class are clear. In Class-0, all the three tables represent no-preference conditions, while in the `Aligned' class, the
preferences in the individual tables are aligned with that in the aggregated table. If the condition of Class-0 occurs, then it can be said with strong confidence that
the data does not provide any choice. If the condition of `Aligned' class occurs, then a choice can be made without any ambiguity. 
Consider the cases of `Weak paradox'. The inequalities in the Tables $T_1$ and $T_2$ are strict and in the same direction, but equality holds in the Table $T_1+T_2$. 
So, it can be said that the preference which can be seen in the individual tables disappears in the aggregated table. 
The more vexing scenario is the cases of `Paradox'. The inequalities in the two individual tables are strict and in the same direction, while the strict inequality 
in the aggregated table is in the opposite direction. This creates the paradox of reversal, where a choice which is individually better in the sub-population becomes 
worse-off in the aggregated population. Making a proper choice in the face of this paradox is a difficult problem.

Informally, Simpson's paradox can be understood to be the following statement covering the classes `Weak paradox' and `Paradox'. 
\begin{quote}
Both the sub-populations provide a clear preference to one particular choice, but this preference is either reversed or ceases to exist in the entire population. 
\end{quote}

Let us now consider the formal definitions of Simpson's paradox in Blyth (1972), Mittal (1991) and Bandyopadhyay et al. (2011).
\begin{description}
	\item{$\mbox{B72}$:} $(A_1\geq C_1) \wedge (A_2\geq C_2) \wedge (\mu<\nu)$ (Blyth (1972)).
	\item{$\mbox{M91}$:} $((A_1\geq C_1) \wedge (A_2\geq C_2) \wedge (\mu\leq\nu)) \vee ((A_1\leq C_1) \wedge (A_2\leq C_2) \wedge (\mu\geq\nu))$ (Mittal (1991).
	\item{$\mbox{BNGBB11}$:} $(A_1\geq C_1) \wedge (A_2\geq C_2) \wedge (\mu\leq \nu)$ and at least one of the inequalities is strict. (Bandyopadhyay et al. (2011)).
\end{description}
The condition M91 is given in Mittal (1991) in a different, but equivalent form. 
The condition B72 captures the Cases~3,~6,~12, and~15, the condition M91 captures the Cases~2,~3,~5,~6,~11 to~17,~22,~23,~25,~26, and the 
condition BNGBB11 captures the Cases~2,~3,~5,~6,~11,~12, and~15. 

The condition M91 covers both the classes `Weak Paradox' and `Paradox'. 
On the other hand, both B72 and BNGBB11 are incomplete, since neither of them capture Case~25 which is a clear case of the paradox. Further, B72 also does not cover
Cases~2 and~26, where a clear preference in the sub-populations ceases to exist in the aggregated population. 
The conditions B72 and BNGBB11 can be expanded. Consider the following conditions.
\begin{description}
	\item {$\mbox{B72}^{\prime}$:} $(A_1\leq C_1) \wedge (A_2\leq C_2) \wedge (\mu>\nu)$. 
	\item {$\mbox{BNGBB11}^{\prime}$:} $(A_1\leq C_1) \wedge (A_2\leq C_2) \wedge (\nu\leq \mu)$ and at least one of the inequalities is strict.
\end{description}
Then $\mbox{Exp-B72}: \mbox{B72}\vee \mbox{B72}^{\prime}$ is the expanded form of Blyth's condition and 
$\mbox{Exp-BNGBB11}: \mbox{BNGBB11}\vee \mbox{BNGBB}^{\prime}$ is the expanded form of Bandyopadhyay et al.'s condition. We note that Exp-B72 still does not
cover Cases~2 and~26 and hence is still incomplete. 

In B72, Exp-B72, BNGBB11, Exp-BNGBB11 and M91, the conditions on the sub-populations (i.e., in tables $T_1$ and $T_2$) include the equality condition. As a consequence,
certain cases get included in the definition which would not be part of the informal description of Simpson's paradox mentioned above. For example, Case~15, 
i.e., $A_1=C_1$, $A_2=C_2$ and $\mu<\nu$ is covered by both B72 and BNGBB11. This case states that both the sub-populations indicate no preference while a preference
is indicated by the whole population. While this may be considered to be surprising, it is certainly not a case where both the sub-populations show a clear trend which 
is reversed or ceases to exist in the entire population. 

The condition M91 has equality conditions on both the sub-populations as well as the aggregated population. This results in M91 covering the Case~14, where equality 
holds for both the sub-populations as well as the aggregate population. This clearly cannot be considered to be surprising. 

Based on the above discussion, we formulate the following definition to exactly cover the informal statement of Simpson's paradox.
\begin{definition}\label{def-SP}
	Given two $2\times 2$ contingency tables $T_1$ and $T_2$, Simpson's paradox (SP) is said to occur for $T_1$ and $T_2$ if 
	the condition $\sym{SP}$ holds where $\sym{SP}=\sym{SP}_1\vee \sym{SP}_2$, with
	\begin{description}
		\item{$\sym{SP}_1$:} $(A_1>C_1) \wedge (A_2>C_2) \wedge \neg(\mu>\nu)$,
		\item{$\sym{SP}_2$:} $(A_1<C_1) \wedge (A_2<C_2) \wedge \neg(\mu<\nu)$.
	\end{description}
\end{definition}
The condition $\sym{SP}_1$ captures Cases~2 and~3 while $\sym{SP}_2$ captures Cases~25 and~26 of Figure~\ref{fig-reln}. 
Suppose $\sym{SP}_2$ holds. Interchanging the first and second rows of the tables $T_1,T_2$ and $T_3$, we see that $\sym{SP}_1$ holds. So, an interchange of the rows of 
the tables converts $\sym{SP}_2$ to $\sym{SP}_1$. Further, both $\sym{SP}_1$ and $\sym{SP}_2$ cannot simultaneously hold. 

Since Exp-B72 does not cover the cases of `Weak paradox', we make further comparison only with M91 and Exp-BNGBB11. Both our definition and Exp-BNGBB11 do not cover the 
classes labelled `Aligned' and the Classes-0,~1,~2 and~5, i.e., cases in these classes are not considered to be examples of SP by either our definition or under Exp-BNGBB11. 
The difference between Exp-BNGBB11 and our definition is that Exp-BNGBB11 includes Classes~3,~4 and~6, while our definition does not. The difference between M91 and
our definition is that M91 includes Classes~0,~3,~4 and~6, while our definition does not. As mentioned earlier, Class~0 is clearly not paradoxical. Next we contrast some of the
classes which are covered by M91 and Exp-BNGBB11 with some of those which are not covered by these conditions.
\begin{enumerate}
\item Contrast Class-4 with Class-1. The similarity
between these two classes is that the inequality in the aggregated population is strict and is in opposition to the strict inequality in one of the sub-populations; 
the dissimilarity is that in Class-1, the inequality in the other sub-population is strict and in reverse direction to that in the aggregated population, while in 
Class-4, equality holds in the other sub-population. This viewpoint suggests that if Class-1 is not viewed as a case of SP, then neither should Class-4 be. 
\item In Class-3, there is a strict inequality in exactly one of the sub-populations, while equality holds in the other as well as the aggregated population. So, a clear 
preference indicated by exactly one of the sub-populations ceases to exist in the aggregated population. Again, we contrast with Class-1, where a clear preference indicated 
by one of the sub-populations is reversed in the aggregated population. So, if Class-1 is not considered to be a case of SP, neither should Class-3 be. 
\item In Class-6, both the sub-populations indicate a no-preference condition, while the aggregated population shows a clear preference. Depending upon the context, this
may or may not be surprising. On the other hand, Class-6 is clearly not covered by the informal statement describing SP.
\end{enumerate}

From Theorem~\ref{thm-possible}, it follows that all of the cases in Table~\ref{fig-reln} are logically possible. Though possible, not all of them are equally likely. 
In particular, the equality conditions in $T_1$, $T_2$ and $T_1+T_2$ are perhaps less likely and hence quite unlikely to occur in practice. At this point though, 
we do not have a formal proof that the equality conditions for $T_1$, $T_2$ and $T_1+T_2$ are less likely to occur compared to the other cases. 
If we indeed assume that equality conditions do not occur in any of these tables, then the definition of SP that we have formulated is equivalent to Exp-B72, M91 and Exp-BNGBB11. 

We have considered SP from the viewpoint of aggregating two contingency tables. The concept easily extends to more than two tables. Our definition of SP has formalised the
notion of cessation or reversal of relation upon aggregation. Interestingly, a scenario of Simpson's paradox has been considered to arise even without cessation or reversal.
Rather, it is based on the magnitudes of the relevant probabilities. For example, if $A_1\gg C_1$ and $A_2\gg C_2$, and $\mu>\nu$, but not $\mu\gg \nu$, then the
strength of the relation in the two sub-populations has been reduced upon aggregation. Such an effect has been considered to be paradoxical. See for
example Page~376 of Pagano and Gauvreau (1993). A formal definition of SP which captures the magnitude effect can be developed and the consequent results analysed. 
We leave such work for the future.

\section{Necessary and/or Sufficient Conditions for Simpson's Paradox \label{sec-condn} }
The consideration for investigating necessary and sufficient conditions for the paradox is inextricably tied to be able to provide the definition of SP,
as giving a definition is typically a philosophical activity.
The additional consideration for it is to consider whether those conditions require any causal intuition to understand 
the paradox as Pearl thinks that causal intuition is the key to unravel its paradoxicality.

A recurrent theme of the research on SP has been to identify either sufficient conditions for SP not to occur, or to identify necessary conditions for SP (so that
taking contrapositive provides sufficient conditions for SP not to occur). Such results have appeared in Lindley and Novick (1981), Mittal (1991)
and Bandyopadhyay et al. (2011). The present section provides a detailed investigation of the necessary conditions for SP and sufficient conditions for SP not to
occur. In the process, we provide new conditions and also simplify and improve upon the proofs of some of the existing conditions. 

Theorem~\ref{thm-SP-nec-cond} below provides necessary conditions for SP to hold and consequently, the negations of these conditions constitute sufficient
conditions for SP not to hold.
\begin{theorem} \label{thm-SP-nec-cond}
	Suppose $\sym{SP}$ occurs for $T_1$ and $T_2$. Then the following conditions hold.
\begin{enumerate}
	\item $A_1\neq A_2$.
	\item $C_1\neq C_2$.
	\item $\alpha_1\neq \gamma_1$ or $\alpha_2\neq \gamma_2$.
\end{enumerate}
\end{theorem}
We note that in TH1 and TH2, Bandyopadhyay et al. (2011) had proved that the Condition BNGBB11 implies $A_1\neq A_2$ and $C_1\neq C_2$. 
The proofs that we obtain for the first two points of Theorem~\ref{thm-SP-nec-cond} are simpler than the proofs of TH1 and TH2 in Bandyopadhyay et al. (2011). 




A concept that turns out to be useful in our study of SP is the well known Odds Ratio which is defined as follows.

\begin{definition}[Odds Ratio]\label{def-OR}
	The Odds (or, Cross-Product) Ratio $\kappa(T)$ for a $2\times 2$ contingency table as in~\eqref{eqn-cont-table} is defined to be $\kappa(T)=(ad)/(bc)$.
\end{definition}

The following result provides a characterisation of Simpson's paradox in terms of the Odds Ratio.
\begin{theorem}\label{thm-SP-basic}
Let $T_1$ and $T_2$ be two $2\times 2$ contingency tables. Simpson's paradox for $T_1$ and $T_2$ is logically equivalent to the following condition.
\begin{eqnarray}\label{eqn-SP-OR}
	& \left( \left(\kappa(T_1) > 1\right) \wedge \left(\kappa(T_2) > 1\right) \wedge \neg\left(\kappa(T_1+T_2) > 1\right) \right) \nonumber \\
	& \vee \nonumber \\
	& \left( \left(\kappa(T_1) < 1\right) \wedge \left(\kappa(T_2) < 1\right) \wedge \neg\left(\kappa(T_1+T_2) < 1\right)\right) .
\end{eqnarray}
\end{theorem}

Based on Theorem~\ref{thm-SP-basic}, we state the following result which provides a characterisation of the condition under which SP does not hold.
\begin{theorem} \label{thm-SP-char}
Let $T_1$ and $T_2$ be two $2\times 2$ contingency tables. Then Simpson's paradox does not hold for $T_1$ and $T_2$ if and only if
\begin{eqnarray}
	\left(\kappa(T_1) > 1\right) \wedge \left(\kappa(T_2) > 1\right) & \mbox{implies} & \kappa(T_1+T_2) > 1, \label{eqn-notSP1} \\
	& \mbox{ and } & \nonumber \\
	\left(\kappa(T_1) < 1\right) \wedge \left(\kappa(T_2) < 1\right) & \mbox{implies} & \kappa(T_1+T_2) < 1 \label{eqn-notSP2}.
\end{eqnarray}
\end{theorem}

\paragraph{Homogeneous subpopulations:} 
Mittal (1991) starts the paper with the following sentence: ``Homogeneity of subpopulations is a very hard concept to define.'' She goes on to explain that
inspite of the difficulty, the concept is required in statistical studies and that the appropriate concept of homogeneity depends on the context. The paper defines
a conception of homogeneity of subpopulations which is shown to be a sufficient condition for SP not to occur. 
We explain the definition of homogeneity introduced in Mittal (1991) with respect to the $2\times 2$ contingency tables in Definition~\ref{defn-cont-tab}.
\begin{definition}[Mittal (1991)]\label{defn-homogeneity}
Let $(T_1,T_2)$ be as defined in Definition~\ref{defn-cont-tab}. Then $(T_1,T_2)$ is said to be homogeneous if one of the following 
conditions hold.
\begin{eqnarray*}
	\max(a_1/b_1,a_2/b_2) < \min(c_1/d_1,c_2/d_2), & &  \max(c_1/d_1,c_2/d_2) < \min(a_1/b_1,a_2/b_2) \\
	\max(a_1/c_1,a_2/c_2) < \min(b_1/d_1,b_2/d_2), & &  \max(b_1/d_1,b_2/d_2) < \min(a_1/c_2,a_2/c_2).
\end{eqnarray*}

\end{definition}
Mittal (1991) uses $\leq$ in the above definition to correspond to her definition of SP, i.e., M91. We have changed this to $<$, to correspond to our definition
of SP as stated in Definition~\ref{def-SP}. The reader may refer to Section~\ref{sec-def-SP} for a discussion on the relation between M91 and Definition~\ref{def-SP}.

The following result is a reformulation of the result from Mittal (1991) showing that the notion of homogeneity captured in Definition~\ref{defn-homogeneity}
is sufficient to rule out SP. The proof provided in Mittal (1991) is geometric. In the appendix, we provide a simple proof based on the Odds Ratio characterisation of SP.
\begin{theorem} \label{thm-suff-m91}
	Let $T_1$ and $T_2$ be two $2\times 2$ contingency tables satisfying Definition~\ref{defn-homogeneity}. Then Simpson's paradox does not hold
	for $T_1$ and $T_2$.
\end{theorem}

Good and Mittal (1987) define a measure of assocation $\alpha$ for the table $T$ as a function $\alpha(T)$ of $T$ satisfying certain desirable properties.
We do not provide the exact properties as these are not directly relevant for the present work. We only note that they showed that the Odds Ratio 
satisfies these properties and is consequently a measure of association according to their definition of the concept.
Given a measure of association $\alpha$, Good and Mittal (1987) define a notion of homogeneity. The two subpopulations
$T_1$ and $T_2$ are homogeneous with respect to $\alpha$ if
$\alpha(T_1) = \alpha(T_2) = \alpha(T_1+T_2)$.
Applying the definition of homogeneity used by Good and Mittal (1987) to the Odds Ratio $\kappa$, the notion of Odds Ratio homogeneity is obtained to be the following.
\begin{definition}[Odds Ratio Homogeneity (ORH)]\label{def-ORH}
The two subpopulations $T_1$ and $T_2$ are homogeneous with respect to the Odds Ratio $\kappa$ if 
\begin{eqnarray}\label{eqn-ORH}
\kappa(T_1) = \kappa(T_2) & = & \kappa(T_1+T_2).
\end{eqnarray}
\end{definition}

We introduce a weaker form of ORH.  
\begin{definition}[Weak Odds Ratio Homogeneity (WORH)] \label{def-WORH}
The two subpopulations $T_1$ and $T_2$ are weakly homogeneous with respect to the Odds Ratio $\kappa$ if 
\begin{eqnarray}\label{eqn-WORH}
\kappa(T_1) = \kappa(T_1+T_2) & \mbox{ or } & \kappa(T_2) = \kappa(T_1+T_2).
\end{eqnarray}
\end{definition}

The next result provides a sufficient condition in terms of WORH under which Simpson's paradox is guaranteed not to hold.
\begin{theorem}\label{thm-WORH-notSP}
Let $T_1$ and $T_2$ be two $2\times 2$ contingency tables. If WORH holds for $T_1$ and $T_2$, then Simpson's paradox
does not hold for $T_1$ and $T_2$.
\end{theorem}

The notion of positive association of true/false valued random variables $X$ and $Y$ is defined in the following manner (Lindley and Novick (1981), Mittal (1991)).

\begin{definition}[Positive Association]\label{def-pos-association}
	The random variables $X$ and $Y$ are said to be positively associated if $\Pr[Y|X] > \Pr[Y|\overline{X}]$. 
\end{definition}
In other words, the conditional probability of $Y$ being true given that $X$ is true is greater than the conditional probability of $Y$ being true given that $X$
is false. So, if we consider two universes (or, sample spaces), one where $X$ is true and the other where $X$ is false, then the probability that $Y$ is true is
greater in the first universe than in the second. 

If $X$ and $Y$ are positively associated, this is denoted as $X\sim Y$. 
Note that as stated, the definition of positive association is asymmetric with respect to the variables $X$ and $Y$, i.e., it is not immediately clear
that $\Pr[Y|X] > \Pr[Y|\overline{X}]$ is equivalent to $\Pr[X|Y] > \Pr[X|\overline{Y}]$.
So, {\em a priori} it is not clear that the 
relation $\sim$ is symmetric\footnote{Mittal (1991) states that if $X\sim Y$ implies $Y\sim X$, then $\sim$ is reflexive; a relation satisfying the stated condition 
is usually called symmetric; a relation $\sim$ is said to be reflexive if $X\sim X$ for all $X$.}.
The next result characterises the notion of positive association in terms of Odds Ratio.

\begin{theorem}\label{thm-odd-ratio-pos-assoc}
	The random variables $X$ and $Y$ are positively associated if and only if $\kappa(T)>1$. Similarly, $X$ and $\overline{Y}$ are positively associated
	if and only if $\kappa(T)<1$. 
\end{theorem}
The proof of the above result for positive association of $X$ and $Y$ shows that $\Pr[Y|X] > \Pr[Y|\overline{X}]$ holds if and only if $\kappa(T)>1$. Similarly, one may show that
$\Pr[X|Y] > \Pr[X|\overline{Y}]$ holds if and only if $\kappa(T)>1$. So, even though the definition of positive association is asymmetric in terms of $X$ and $Y$, the
relation $X\sim Y$ is actually symmetric. Seen in terms of Odds Ratio, the symmetric property of $\sim$ becomes immediate. Showing this property required a somewhat more
complicated proof in Mittal (1991).
Note that the condition $X$ is positively associated with $\overline{Y}$ is not the same as the negation of the condition that $X$ is positively associated with $Y$ since the
possibility $\kappa(T)=1$ is missing from the characterisations of both the conditions.

Given a pair of contingency tables $T_1$ and $T_2$ as in Definition~\ref{defn-cont-tab} and the associated random variables as described in the discussion following
the definition, if $X$ is positively associated with $Y$
in table $T_1$ (resp. $T_2$), we denote this as $(X\sim Y)|M$ (resp. $(X\sim Y)|\overline{M}$). The notation $(X\sim Y)|M$ can be understood as saying that
$X$ is positively associated with $Y$ with respect to the sub-population represented by $T_1$, and similarly the notation $(X\sim Y)|\overline{M}$ can be understood as saying that
$X$ is positively associated with $Y$ with respect to the sub-population represented by $T_2$. In particular, $X$ and $Y$ being positively associated with respect to
one of the sub-populations represented by $T_1$ or $T_2$ does not necessarily imply that $X$ and $Y$ are also positively associated with the other sub-population or, 
with the overall population. 
We write $X\sim Y$ if the positive association condition holds for the merged table $T_1+T_2$.

Using Theorem~\ref{thm-odd-ratio-pos-assoc}, $(X\sim Y)|M$ (resp. $(X\sim Y)|\overline{M}$) holds if and only if $\kappa(T_1)>1$ (resp. $\kappa(T_2)>1$); and
$X\sim Y$ if and only if $\kappa(T_1+T_2)>1$. 

Theorem~\ref{thm-SP-basic} provides a condition in terms of Odds Ratio which is equivalent to the definition of Simpson's paradox. The connection of positive
association to Odds Ratio can be used to translate this condition in terms of positive association of the random variables $X$ and $Y$. Using
Theorems~\ref{thm-SP-basic} and~\ref{thm-odd-ratio-pos-assoc}, SP holds for $(T_1,T_2)$ if and only if either 
($(X\sim Y)|M$ and $(X\sim Y)|\overline{M}$ and $\neg(X\sim Y)$), or ($(X\sim \overline{Y})|M$ and $(X\sim \overline{Y})|\overline{M}$ and $\neg(X\sim \overline{Y})$).
In other words, we have the following result.
\begin{theorem}\label{thm-SP-positive-association}
	SP holds if and only if {\em either} $X$ is positively associated with $Y$ in both the sub-populations and not so in the aggregate
population, {\em or} $X$ is positively associated with $\overline{Y}$ in both the sub-populations and not so in the aggregate population.
\end{theorem}

The unconditional probabilities of the events $M\wedge Y$, $M\wedge \overline{Y}$, $\overline{M}\wedge Y$, $\overline{M}\wedge \overline{Y}$ and
$M\wedge X$, $M\wedge \overline{X}$, $\overline{M}\wedge X$, $\overline{M}\wedge \overline{X}$ are obtained respectively from the following $2\times 2$ contingency
tables $S_1$ and $S_2$ given below. 
Note that $S_1$ and $S_2$ are obtained from the tables $T_1$ and $T_2$ in Definition~\ref{defn-cont-tab} by aggregating in a different manner.
In both $S_1$ and $S_2$, the variables $M$ and $\overline{M}$ are associated with the first and second rows respectively. In $S_1$, $Y$ and $\overline{Y}$ are
associated with the first and second columns respectively, while in $S_2$, $X$ and $\overline{X}$ are associated with the first and second columns respectively.
\begin{center}
\begin{tabular}{cc}
\begin{minipage}{150pt}
\begin{eqnarray*}
S_1 & = & 
\begin{array}{|c|c|}
	\hline
	a_1+c_1 & b_1+d_1 \\ \hline
	a_2+c_2 & b_2+d_2 \\ \hline
\end{array}\ ,
\end{eqnarray*}
\end{minipage}
\begin{minipage}{150pt}
\begin{eqnarray*}
S_2 & = & 
\begin{array}{|c|c|}
\hline
a_1+b_1 & c_1+d_1 \\ \hline
a_2+b_2 & c_2+d_2 \\ \hline
\end{array}\ .
\end{eqnarray*}
\end{minipage}
\end{tabular}
\end{center}
Using Theorem~\ref{thm-odd-ratio-pos-assoc}, $M\sim Y$ (resp. $M\sim X$) if and only if $\kappa(S_1)>1$ (resp. $\kappa(S_2)>1$); and
$M\sim \overline{Y}$ (resp. $M\sim \overline{X}$) if and only if $\kappa(S_1)<1$ (resp. $\kappa(S_2)<1$).

Lindley and Novick (1981) had proposed a necessary condition for SP. The result of Lindley and Novick (1981) was later expanded by Mittal (1991). 
Mittal (1991) had pointed out a few problematic issues in the arguments used by Lindley and Novick (1981)
and had provided a new proof. The proof given in Mittal (1991) itself is essentially a sketch and omits some important details. Also,
there are some inaccuracies in the handling of the inequalities. 
To resolve these issues and for the sake of completeness, we provide a detailed proof of the result in the appendix.

\begin{theorem}[Lindley and Novick (1981), Mittal (1991)] \label{thm-nec-cond}
Let $T_1$ and $T_2$ be a pair of contingency tables and random variables $X$, $Y$ and $M$ are associated with $T_1$ and $T_2$ as stated above. If
	$(X\sim Y)|M$ and $(X\sim Y)|\overline{M}$ and SP occurs for $(T_1,T_2)$, then either (a) $M\sim X$ and $M\sim Y$; or (b) $M\sim \overline{X}$ and $M\sim \overline{Y}$.
	Further, the result also holds with $Y$ replaced by $\overline{Y}$.
\end{theorem}

We stated several theorems providing necessary and/or sufficient conditions for SP. Pearl contends that the only way to untangle the paradox is to fall back on causal intuitions. 
The theorems about the paradox presented in this section show that his contention is not necessarily true. 
The paradox, or at least some versions of it, can be understood without requiring causal intuitions (see Section~\ref{sec-ex}  for more on this theme).

\section{Simpson's Paradox and Progressive Accumulation of Data \label{sec-data} }
Statistics has been viewed as a method of inductive inference. We mention some examples from the literature where such a connection has been pointed out by leading
statisticians. Mahalanobis (1950) wrote: ``The importance of statistics in the field of science is due to its supplying the general method for inductive inference.''
Fisher (1955) had written about the relation between statistical methods and scientific induction. The preface of the book by Rao (1997) 
suggests that statistics provides a method of codifying inductive reasoning. Mayo and Cox (2006) explore the role of frequentist statistics as a theory of inductive
inference. 

It is well known that inductive inference is non-monotonic (see for example Levi (2005)). In the present context, the pair of $2\times 2$ contingency tables
can be considered to be data. An inference based on such a pair of tables is essentially an inductive inference. Due to the non-monotonic nature of inductive
inference, an inference drawn from a particular pair of tables may get invalidated as more data are obtained. One may consider a scenario where data are progressively
accumulated over time. Inferences are made at certain time points based on the current data that are available. This leads to the theoretical possibility that
the data indicates SP at a certain point of time which disappears as more data accumulates and then reappears with further accumulation of data. We present a
result to show that such a theoretical possibility can indeed arise and the appearance and disappearance of SP can continue ad infinitum. 

The importance of the role of inductive inference in both philosophy and statistics is well appreciated. By noting the inductive nature of SP, we will show that the 
most significant question about the paradox, that is, whether to aggregate or not, when confronted with a SP type situation, often rests on the appearance and 
disappearance of SP, as data accumulate. The decision regarding the what-to-do question in a dynamic SP type situation is yet to be well-appreciated in the literature.

Let $\mathcal{T}=(T_1^{(k)},T_2^{(k)})_{k\geq 0}$ be a sequence of pairs of $2\times 2$ contingency tables, where
\begin{center}
\begin{tabular}{cc}
\begin{minipage}{120pt}
\begin{eqnarray*}
	T_1^{(k)} & = & 
\begin{array}{|c|c|}
\hline
	a_1^{(k)} & b_1^{(k)} \\ \hline
	c_1^{(k)} & d_1^{(k)} \\ \hline
\end{array}\ ,
\end{eqnarray*}
\end{minipage}
\begin{minipage}{120pt}
\begin{eqnarray*}
	T_2^{(k)} & = & 
\begin{array}{|c|c|}
\hline
	a_2^{(k)} & b_2^{(k)} \\ \hline
	c_2^{(k)} & d_2^{(k)} \\ \hline
\end{array}\ .
\end{eqnarray*}
\end{minipage}
\end{tabular}
\end{center}
We define $\mathcal{T}$ to be monotonic if for all $k\geq 0$,
\begin{eqnarray*}
	& a_1^{(k+1)} \geq a_1^{(k)}, \quad b_1^{(k+1)} \geq b_1^{(k)}, \quad c_1^{(k+1)} \geq c_1^{(k)}, \quad d_1^{(k+1)} \geq d_1^{(k)}, \mbox{ and} \\
	& a_2^{(k+1)} \geq a_2^{(k)}, \quad b_2^{(k+1)} \geq b_2^{(k)}, \quad c_2^{(k+1)} \geq c_2^{(k)}, \quad d_2^{(k+1)} \geq d_2^{(k)}. 
\end{eqnarray*}
In other words, the pair of contingency tables $(T_1^{(k+1)},T_2^{(k+1)})$ is obtained from the pair $(T_1^{(k)},T_2^{(k)})$ by accumulating additional data. From 
an empirical data perspective, some initial data collection leads to the pair of tables $(T_1^{(0)},T_2^{(0)})$. Obtaining additional data over and above that which has 
been already accounted for in the pair $(T_1^{(0)},T_2^{(0)})$ leads to the pair of tables $(T_1^{(1)},T_2^{(1)})$, further data successively leads to the 
pairs of tables $(T_1^{(2)},T_2^{(2)})$, $(T_1^{(3)},T_2^{(3)})$, and so on.

Consider Figure~\ref{fig-reln}. Case~3 represents $(>,>,<)$, i.e., the `$>$' relation holds for both $T_1$ and $T_2$, but `$<$' holds for $T_1+T_2$. Similarly,
Case~25 represents $(<,<,>)$. Both Cases~3 and~25 are instances of SP. For the sake of convenience, let us call Case~3 to be positive SP and Case~25 to be
negative SP. Similary, Case~1 represents $(>,>,>)$ and we will call this to be positive alignment and Case~27 represents $(<,<,<)$ which we call to be negative
alignment. 

\begin{theorem}\label{thm-togglings}
	There exists a monotonic sequence $\mathcal{T}=(T_1^{(k)},T_2^{(k)})_{k\geq 0}$ such that the following holds.
\begin{enumerate}
	\item $(T_1^{(0)},T_2^{(0)})$, $(T_1^{(4)},T_2^{(4)})$, $(T_1^{(8)},T_2^{(8)}), \ldots$ are instances of positive SP.
	\item $(T_1^{(1)},T_2^{(1)})$, $(T_1^{(5)},T_2^{(5)})$, $(T_1^{(9)},T_2^{(9)}), \ldots$ are instances of positive alignment.
	\item $(T_1^{(2)},T_2^{(2)})$, $(T_1^{(6)},T_2^{(6)})$, $(T_1^{(10)},T_2^{(10)}), \ldots$ are instances of negative SP.
	\item $(T_1^{(3)},T_2^{(3)})$, $(T_1^{(7)},T_2^{(7)})$, $(T_1^{(11)},T_2^{(11)}), \ldots$ are instances of negative alignment.
\end{enumerate}
\end{theorem}
Theorem~\ref{thm-togglings} shows that it is possible to have progressively accumulated data where positive SP and negative SP are interleaved with positive and negative alignments
and such togglings repeat ad infinitum. Consider $(T_1^{(0)},T_2^{(0)})$. Since this is an instance of positive SP, an empirical scientist may invoke some 
resolution mechanism and draw a conclusion. If the conclusion is to accept the relation suggested by the sub-populations, then further data as in 
$(T_1^{(1)},T_2^{(1)})$ will confirm such a conclusion, while if the conclusion is to accept the relation suggested by the aggregate population, then 
$(T_1^{(1)},T_2^{(1)})$ will invalidate the conclusion. Further data will exhibit negative SP followed by negative alignment and then again positive SP and so on.
This suggests that any explanation of Simpson's paradox which aims to be comprehensive should keep in focus the inductive nature of inference from data.

Progressive accumulation of data is a practical reality. The ongoing Covid-19 pandemic is a very relevant example. As the disease progresses over time, increasingly more data
becomes available. SP may be observed in the data at a particular point of time. For example, Kugelgen, Gresele and Scholkopf (2021) analysed data up to a certain point of 
time to demonstrate that while the overall fatality rate was higher in Italy than in China, within every age group, the fatality rate was lower in Italy. So, a clear
case of SP occurs in the data that was considered by the authors at the time of the analysis. It is not clear that with further accumulation of data, the occurrence
of SP will persist. Answering the what-to-do question in the context of dynamic on/off toggling of SP with progressive accumulation of data does not seem to have received 
sufficient serious attention in the literature.

\section{Deconstructing Examples of Simpson's Paradox \label{sec-ex} }

The literature provides several examples of SP. In this section, we discuss some of these examples. There are two distinct aspects to the examples of SP that appear
in the literature. The first aspect is numerical which pertains to the actual data available in a pair of $2\times 2$ contingency tables, while the second aspect is the 
context in which the contingency tables are discussed. In our discussion, we separate the numerical aspect from the contextual aspect. This leads to the possibility
of combining the numerical aspect of an example with the context of another. As a result, it is possible to argue that a contextual explanation that is meaningful
with respect to a particular numerical example does not remain so with respect to another numerical example. After discussing the examples, we illustrate this issue
by combining the numerical aspect of one example with the context explanation of another. The point of such an exercise is to highlight that any explanation of 
SP must combine both the data and the context underlining the fundamental empirical nature of SP.
Providing many examples of the same paradox is to argue that because of numerical and contextual features of the paradox, a unified treatment of them under one umbrella 
may not be able to do justice to the complexity and various dimensions of SP.

Pearl (2014) had commented that Lindley and Novick (1981) had raised ``Simpson's paradox to new heights'' by providing two contexts for the same data. While Lindley and Novick's 
work was indeed influential, two separate contexts for the same data appear in the much earlier paper by Simpson (1951). 

We describe the context in terms of the true/false valued random variables $X,Y$ and $M$ identified in Section~\ref{sec-def-SP}. 
To refresh memory, we mention that $M$ and $\overline{M}$ correspond to the tables $T_1$ and $T_2$ respectively; $X$ and $\overline{X}$ correspond to the first and
the second rows respectively; and, $Y$ and $\overline{Y}$ correspond to the first and second columns respectively.

\paragraph{Example~1 (Simpson 1951):} 
\begin{center}
\begin{tabular}{cc}
\begin{minipage}{100pt}
\begin{eqnarray*}
T_1 & = & 
\begin{array}{|c|c|}
\hline
4 & 3 \\ \hline
8 & 5 \\ \hline
\end{array}\ ,
\end{eqnarray*}
\end{minipage}
\begin{minipage}{100pt}
\begin{eqnarray*}
T_2 & = & 
\begin{array}{|c|c|}
\hline
2 & 3 \\ \hline
12 & 15 \\ \hline
\end{array}\ ,
\end{eqnarray*}
\end{minipage}
\begin{minipage}{100pt}
\begin{eqnarray*}
T_1+T_2 & = &
\begin{array}{|c|c|}
\hline
6 & 6 \\ \hline
20 & 20 \\ \hline
\end{array}\ .
\end{eqnarray*}
\end{minipage}
\end{tabular}
\end{center}
For this example, we have $A_1<C_1$, $A_2<C_2$ and $\mu=\nu$. So, the association which is visible in the sub-populations ceases (as opposed to being reversed)
upon aggregation.

\paragraph{Context: Cards (Simpson 1951).} The context concerns a deck of cards. 
$M$ and $\overline{M}$ correspond to dirty and clean cards respectively;
$X$ and $\overline{X}$ correspond to court (King, Queen, Knave) cards and plain cards respectively; $Y$ and $\overline{Y}$ correspond to red and black cards respectively. 
The goal is to examine whether the proportion of court cards is associated with colour. It would appear to be so from the individual tables $T_1$ and $T_2$, but
such association disappears from the aggregated table. As Simpson noted, in this context, it is the aggregated table which makes sense. 

\paragraph{Context: Treatment (Simpson 1951).} The context concerns whether patients are alive or dead after receiving treatment or not. 
$M$ and $\overline{M}$ correspond to male and female patients respectively;
$X$ and $\overline{X}$ correspond to untreated and treated patients respectively; $Y$ and $\overline{Y}$ correspond to being alive and dead respectively. 
So, a positive association between treatment and survival is seen separately for male and female population groups, but no such association is observed
from the aggregated data. Simpson raises the question of ``sensible'' interpretation noting that a treatment which is beneficial to both males and females cannot
be rejected as being useless for the race. The wording of the question suggests that Simpson considers the inference from the disaggregated tables to be more sensible, but
does not actually state this.

\paragraph{Example~2 (Blyth 1971):} 
\begin{center}
\begin{tabular}{cc}
\begin{minipage}{120pt}
\begin{eqnarray*}
T_1 & = & 
\begin{array}{|c|c|}
\hline
1000 & 9000 \\ \hline
50 & 950 \\ \hline
\end{array}\ ,
\end{eqnarray*}
\end{minipage}
\begin{minipage}{120pt}
\begin{eqnarray*}
T_2 & = & 
\begin{array}{|c|c|}
\hline
95 & 5 \\ \hline
5000 & 5000 \\ \hline
\end{array}\ ,
\end{eqnarray*}
\end{minipage}
\begin{minipage}{120pt}
\begin{eqnarray*}
T_1+T_2 & = &
\begin{array}{|c|c|}
\hline
1095 & 9005 \\ \hline
5050 & 5950 \\ \hline
\end{array}\ .
\end{eqnarray*}
\end{minipage}
\end{tabular}
\end{center}
For this example, we have $A_1>C_1$, $A_2>C_2$ and $\mu<\nu$. So, the association which is visible in the sub-populations is reversed upon aggregation.

\paragraph{Context: Treatment (Blyth 1971).} The context concerns whether local or outstation patients are treated. 
$M$ and $\overline{M}$ correspond to local and outstation patients respectively;
$X$ and $\overline{X}$ correspond to a new treatment and the standard treatment respectively; 
$Y$ and $\overline{Y}$ correspond to being alive and dead respectively.
Blyth explains the reversal by suggesting that local patients are much less likely to recover and the treatment was given mostly to local patients. A treatment will
show a poor recovery rate if it is given mostly to the most seriously ill patients. 

\paragraph{Example~3 (Gardner 1976):} 
\begin{center}
\begin{tabular}{cc}
\begin{minipage}{100pt}
\begin{eqnarray*}
T_1 & = & 
\begin{array}{|c|c|}
\hline
5 & 6 \\ \hline
3 & 4 \\ \hline
\end{array}\ ,
\end{eqnarray*}
\end{minipage}
\begin{minipage}{100pt}
\begin{eqnarray*}
T_2 & = & 
\begin{array}{|c|c|}
\hline
6 & 3 \\ \hline
9 & 5 \\ \hline
\end{array}\ ,
\end{eqnarray*}
\end{minipage}
\begin{minipage}{100pt}
\begin{eqnarray*}
T_1+T_2 & = &
\begin{array}{|c|c|}
\hline
11 & 9 \\ \hline
12 & 9 \\ \hline
\end{array}\ .
\end{eqnarray*}
\end{minipage}
\end{tabular}
\end{center}
For this example, we have $A_1>C_1$, $A_2>C_2$ and $\mu<\nu$. So, the association which is visible in the sub-populations is reversed upon aggregation.

\paragraph{Context: Poker Chips (Gardner 1976).} The context concerns different coloured poker chips placed in different coloured hats. 
$M$ and $\overline{M}$ correspond to two tables (i.e., pieces of furniture, not tabulation of data);
$X$ and $\overline{X}$ correspond to black and grey hats respectively;
$Y$ and $\overline{Y}$ correspond to being blue and white chips respectively.

\paragraph{Example~4 (Lindley and Novick 1981):} 
\begin{center}
\begin{tabular}{cc}
\begin{minipage}{100pt}
\begin{eqnarray*}
T_1 & = & 
\begin{array}{|c|c|}
\hline
18 & 12 \\ \hline
7 & 3 \\ \hline
\end{array}\ ,
\end{eqnarray*}
\end{minipage}
\begin{minipage}{100pt}
\begin{eqnarray*}
T_2 & = & 
\begin{array}{|c|c|}
\hline
2 & 8 \\ \hline
9 & 21 \\ \hline
\end{array}\ ,
\end{eqnarray*}
\end{minipage}
\begin{minipage}{100pt}
\begin{eqnarray*}
T_1+T_2 & = &
\begin{array}{|c|c|}
\hline
20 & 20 \\ \hline
16 & 24 \\ \hline
\end{array}\ .
\end{eqnarray*}
\end{minipage}
\end{tabular}
\end{center}
For this example, we have $A_1<C_1$, $A_2<C_2$ and $\mu>\nu$. So, the association which is visible in the sub-populations is reversed upon aggregation.

\paragraph{Context: Treatment (Lindley and Novick 1981).} The context concerns treatments to patients.
$M$ and $\overline{M}$ correspond to male and female patients respectively;
$X$ and $\overline{X}$ correspond to being treated and not treated respectively;
$Y$ and $\overline{Y}$ correspond to recovery or not respectively.
In this context, not being treated appears to be better individually for both male and female and patients, but upon aggregation, it appears that treatment is better. The
explanation forwarded by Lindley and Novick is that males have been mostly assigned to the treatment group and females to the control group. They also suggested that
perhaps the doctor distrusted the treatment and was reluctant to give it to females where the recovery rate is much lower. 
For this context, the disaggregated data makes more sense, since it would be better not to recommend the treatment when it is clearly inferior for both males and females.
Lindley and Novick had further suggested that
instead of male/female distinction, the groupings can be considered based on some factor which may be difficult to determine such as a genetic classification. 

\paragraph{Context: Agriculture (Lindley and Novick 1981).} The context concerns coloured varities of plants which are short or tall and have low or high yield.
$M$ and $\overline{M}$ correspond to tall and short plants respectively;
$X$ and $\overline{X}$ correspond to white and black varieties respectively;
$Y$ and $\overline{Y}$ correspond to high and low yield respectively.
For this context, the aggregated data makes more sense, since one would choose the white variety which provides higher yield overall ignoring the individual yield
information provided by the tall and the short varities.

\paragraph{Example~5 (Hand 1994):} 
\begin{center}
\begin{tabular}{cc}
\begin{minipage}{100pt}
\begin{eqnarray*}
T_1 & = & 
\begin{array}{|c|c|}
\hline
255 & 174 \\ \hline
156 & 102 \\ \hline
\end{array}\ ,
\end{eqnarray*}
\end{minipage}
\begin{minipage}{100pt}
\begin{eqnarray*}
T_2 & = & 
\begin{array}{|c|c|}
\hline
88 & 222 \\ \hline
82 & 175 \\ \hline
\end{array}\ ,
\end{eqnarray*}
\end{minipage}
\begin{minipage}{100pt}
\begin{eqnarray*}
T_1+T_2 & = &
\begin{array}{|c|c|}
\hline
343 & 396 \\ \hline
238 & 277 \\ \hline
\end{array}\ .
\end{eqnarray*}
\end{minipage}
\end{tabular}
\end{center}
For this example, we have $A_1<C_1$, $A_2<C_2$ and $\mu>\nu$. So, the association which is visible in the sub-populations is reversed upon aggregation.

\paragraph{Context: Psychiatry (Hand 1994).} The context concerns the proportion of male and female patients of different age groups in psychiatry wards in two
different years.
$M$ and $\overline{M}$ correspond to patients having ages $\leq 65$ and those having ages $>65$ respectively;
$X$ and $\overline{X}$ correspond to numbers of patients in the years 1970 and 1975 respectively;
$Y$ and $\overline{Y}$ correspond to male and female respectively.
Hand contends that both the disaggregated and the aggregated table provide meaningful answers, depending on what question one asks. If one wishes to know whether the
proportion of males increase, then the aggregated table provides the answer, while, if one wishes to know for patients of a given age, whether the proportion of males
increase, then the disaggregated tables provide the answer.

Our deconstruction of the examples of SP into the numerical and the contextual aspect raises interesting possibilities. One may take the numerical aspect
of one example and match with the contextual example of another example. This gives rise to a new example of SP. One may then contrast the explanation of the context
provided for the original numerical part with that of the new numerical part. 
To illustrate this procedure, consider combining the Lindley and Novick (1981) treatment
context with the numerical example of Blyth (1971). One then notes that recovery rate for treated patients are higher than those for untreated patients in both male
and female groups, but overall the recovery rate is lower for treated patients. Should one then go by the individual tables and recommend the treatment? Since the
recovery rate for the aggregated data is lower for treatment, can one simply ignore this finding without further exploration of possible confounding factors? 
In the medical context, it is perhaps better to err on the side of safety and not recommend the treatment. 
In contrast to these doubts, the Lindley-Novick treatment context combined with the Lindley-Novick's numerical example, led to the suggestion of going by disaggregated 
tables and this was considered reasonable. The point here is that an explanation or resolution of SP is for the holistic combination of data-plus-context. Considering an
explanation to be only for the context while ignoring the data may not be satisfying.

\section{The Surprise Element \label{sec-surprise} }
The reversal property in a pair of contingency table is called a paradox, since it appears to be surprising. 
Following the distinction in Section~\ref{sec-ex}, we consider the numerical and the contextual aspect to be separate. There is nothing surprising or unsurprising about the
numerical aspect. To highlight that the numerical aspect is not surprising, we mention the connection of SP to determinants of $2\times 2$ matrices. Consider
the tables $T_1$ and $T_2$ to be $2\times 2$ matrices. The condition $A_1>C_1$, $A_2>C_2$ and $\mu<\nu$ is equivalent to the condition
$\mbox{det}(T_1)>0$, $\mbox{det}(T_2)>0$ and $\mbox{det}(T_1+T_2)<0$, where $\mbox{det}(T)$ denotes the determinant of $T$. (This is easily seen through the connection 
to the Odds Ratio characterisation given in Theorem~\ref{thm-SP-basic}.) There is nothing surprising in the fact that the determinants of $T_1$ and $T_2$ are
both positive while the determinant of $T_1+T_2$ is negative.

In view of the above, the surprise element arises from the context. We provide an explanation of why such a reversal
appears to be surprising. Human intuition is trained to expect {\em uniformity} in various aspects of reasoning. For example, in elementary deductive logic, suppose 
that a statement $P_1$ implies a statement $Q$ and another statement $P_2$ also implies $Q$. The statements $P_1$ and $P_2$ provide support to the inference of $Q$.
Aggregation of the supports will mean $P_1\vee P_2$ and then we still have $P_1\vee P_2$ implies $Q$. More generally, human intuition expects that if some property
holds for parts, then if the parts are put together the same property will still continue to hold for the whole. 

Such uniformity, on the other hand, does not necessarily hold in all possible scenarios. Whenever parts behave differently from the whole, i.e., all the parts display 
a particular behaviour while the whole does not display such a behaviour or, displays a behaviour which is in contrast to the behaviour of the parts, a surprise element is
involved and this gives rise to a paradox. Simpson's paradox is clearly of this kind. A property which holds for the two sub-populations ceases to hold or reverses
when the populations are combined. This can be seen as falling within the broad class of aggregation paradoxes. There are some well known examples of such paradoxes where 
a feature which holds for individuals reverses itself upon aggregation. As examples, we refer to the Ostrogorski paradox (Daudt and Rae (1976)) and the discursive 
dilemma (Kornhauser and Sager (1986), List and Petit (2002)) where choices/judgements made at the individual level are not preserved upon aggregation. Another famous paradox 
involving parts and whole is the Banach-Tarski paradox, which shows that a unit sphere
can be divided into finitely many parts and assembled together to obtain two unit spheres; the paradox being explained by the fact that the individual
parts are defined in a manner that precludes these parts from having a volume, whereas the whole sphere has a well-defined volume.

Lewis wrote, ``[w]hen commonsense delivers a firm and uncontroversial answer about a not-too-far-fetched case, theory had better agree.'' It is hard to say whether 
human intuition and expectations that the uniformity about the whole is combined from its parts is causally motivated. But, at the same time, we cannot brand this intuition as 
a ``far-fetched case.'' Lewis also reminds us that ``[i]f an analysis of causation does not deliver the common-sense answer [then] that is bad trouble.'' 
Gardner's hat example is a clear case which does not involve any causal intuition. It is unclear whether Lewis' dictum should serve as a bench-mark of any account of causation. 
Therefore, the complex nature of SP should neither support nor weaken his dictum.

While explaining Simpson's paradox from the viewpoint of causal graphs, Pearl (2014) comments that 
\begin{quote}
``... it is hard, if not impossible, to explain the surprise part of Simpson's reversal without postulating that human intuition is governed by causal calculus 
	together with a persistent tendency to attribute causal interpretation to statistical associations.'' 
\end{quote}
We consider two objections to the above statement. The first objection is that there is indeed a way to explain the surprise element without assuming that
causal considerations are woven into human intuition.
The second objection is that there are examples of SP in the literature which do not readily admit a causal modelling. The effect of the second objection is to show
the incompleteness of the causal approach in explaining the surprise element of SP while the first objection shows a different method of explaining this surprise element.

Regarding the first objection, we note that as argued above, there is a different and perhaps more fundamental perspective which explains the surprise element of SP, i.e., human 
intuition expects a kind of uniformity where the whole behaves in a manner similar to the parts. Deviations create surprise and lead to other paradoxes as mentioned above.
As for the second objection, in contrast to what Pearl suggests, there are examples of SP where there is no natural interpretation using causal models. 
One such example, is the poker chip context of Gardner (1976) which has been described in Section~\ref{sec-ex}.
In this example, there are no natural causal issues. It is absurd to say that the colours of the hats cause the colours of the chips or vice versa; similarly, it is absurd to
say that the physical tables upon which the hats have been placed cause either the colours of the hats or the colours of the chips. Since a causal diagram is a model, one
may of course, draw any model that one wishes, but there would be no way to determine that one particular model is to be preferred to another. A scenario where all models
are equally (un)likely hardly provides a useful theory. In effect, viewing this example using causal models is like trying to fit a square peg into a round hole. 

Another example which does not admit a causal model has been described by Bandyopadhyay et al. (2011). In this example, there are two bags with each bag containing
marbles which are either big or small and are either red or blue. Though a numerical example is not provided, it is possible to construct an example of a pair of $2\times 2$
contingency tables corresponding to the bags and rows/columns corresponding to size/colour of marbles which exhibit SP. Again, there are no natural causal issues
in this example: bags do not cause size/colour of marbles; size does not cause colour and vice versa. 

The above examples show that viewing SP through the prism of causal models potentially leaves out examples of SP. Note that our explanation of the surprise element arising
from the expectation of uniformity applies to both the above examples.

\section{Causal Structures and SP: Some Questions \label{sec-causal} }
Causal structures have been put forward as the resolution of Simpson's paradox (Pearl 2009). In this section, we put forward questions regarding causal
structures and their relation to SP. Some of these questions are about causal structures themselves and have been raised earlier without being satisfactorily addressed
while the other questions are new and seem difficult to answer. 

Consider a set of random variables $V$. Based on prior knowledge, the variables are ordered by hypothesised cause-effect relations. The variables and the cause-effect relations
can be conveniently described using a directed acyclic graph (DAG), where the nodes are the variables and there is an arrow from one node to another if the first node
is a cause of the second. The joint probability distribution of the random variables in $V$ is denoted by $P(v)$. A causal structure consists of the DAG and the joint
probability distribution. Such a structure was introduced by Pearl (1995, 2009). For a random variable $X$, and a value $x$ in the domain of $X$, the event $X=x$ is well 
defined and one may consider conditional probabilities of the form $\Pr[Y=y|X=x]$, where $Y$ is another random variable in $V$ and $y$ is a value in the domain of $Y$. 
Pearl (2009) defines an important concept called ``do'', where $\mbox{do}(x)$
corresponds to setting $X=x$ everywhere and removing the variable $X$ from the DAG. This alters the structure of the DAG and the joint probability distribution. 
One can consider the probability of the event $Y=y$ in the altered causal structure and this probability is denoted as $\Pr[Y=y|\mbox{do}(x)]$. In Pearl (1995) this
notion was termed intervention. An intervention is not always possible in practice. A remarkable result proved in Pearl (1995) shows that under certain conditions the effect 
of intervention can be determined from the causal structure without actually performing the intervention. 

Pearl (1995, 2009) defines an important concept related to causal graphs, namely the back-door criteria. Given a pair of variables $(X,Y)$, another variable $Z$ 
(or, possibly a set of variables) satisfies the back-door criterion relative to $(X,Y)$, if $Z$ is not a direct descendant of $X$, and $Z$ blocks every path between
$X$ and $Y$ which contains an arrow into $X$. The back-door criterion captures the idea that the variable $Z$ blocks ``spurious'' paths entering $X$. A corresponding
notion of front-door criterion has also been defined. 

While SP was not considered in Pearl (1995), later work by Pearl (2009, 2014) applied 
the causal structure described above and the notion of back-door criterion to the analysis of SP. Certain graph structures are
stated to be incapable of exhibiting SP which in effect constitute sufficient conditions for not SP. Such sufficient conditions are different in nature from the
sufficient conditions for not SP which have been derived in Section~\ref{sec-condn}. The back-door criterion has been put forward by Pearl as a resolution mechanism
for SP. Simple rules are provided for determining whether to accept the conclusion of the individual data, or the aggregated data, or, none. 

The causal modelling and the subsequent analysis of SP is a remarkable achievement. Nonetheless, there are some relevant questions regarding this approach which 
do not seem to have been adequately addressed in the literature. Below we list down some such questions.

\paragraph{Can all examples of SP be modelled using causal structures?} In Section~\ref{sec-surprise}, we have provided two scenarios where there is no natural associated causal
model. This shows that not all examples of SP have a natural interpretation in terms of causal structures.

\paragraph{Can all examples of SP be resolved using causal structures?} Since all examples of SP cannot be modelled using causal structures, it clearly follows that all
examples of SP cannot be resolved using such structures. For examples of SP which do admit causal modelling, there is an implicit claim in Pearl (2014) that a resolution
can be obtained within the causal model. There does not, however, seem to be a proof of such a claim and it is not even clear whether such a claim can be formally stated 
and proved.

\paragraph{Use of deductive methods to resolve an inductive inference.} It has been argued in Section~\ref{sec-data} that statistical inference is in general a 
method of inductive inference. A causal structure on the other hand, assumes that the causal relations are known before hand. This raises the question of how
one arrives at these causal relations? Such relations may be considered to have been provided by domain experts. This, however, deflects the question to another level,
since one could pose the question as to how the domain experts would have arrived at these causal relations. From a scientific point of view, we see no way of justifying
the determination of causal relations except by referring to previous observations/data. But, then one would have to look for justifications of the procedures
for deriving the causal relations from data, leading to a circularity in the argument. An argument along the same line has been given earlier in Freedman (1995) who remarked that
``Validation of causal models remains an unsolved problem'' and also ``... causal laws, ..., are assumed rather than inferred from the data.'' 

At several places, Pearl has critiqued statisticians for avoiding terminology based on causes. For a discussion on use of `correlation' by Pearson and Yule,
we refer to Aldrich (1995). The paper also identifies several aspects of Pearson's thoughts on causality and remarks that Pearson's ``complete position was that
causation {\em is} correlation, or more precisely the limiting case of correlation, {\em except} when the correlation is spurious, when correlation is not
causation.'' (Emphasis as in the original.) It is perhaps not surprising that the field of statistical inference will be careful in talking about causation. After 
all, statistical inference is a form of inductive inference and the skepticism to inferring causes using inductive inference has a long history. 

A causal model incorporates known relationships among the variables. How does one determine whether the known relationships
are {\em all} the relations that hold? Along this line of reasoning, it has been observed by Cox and Wermuth (1995) that the ``back-door criterion requires there
to be no `common cause' ... that is not blocked out by the observed variable'' and that it is ``precisely doubt about such assumptions that makes epidemiologists, 
for example, wisely in our view, so cautious in distinguishing risk factor from causal effects.'' One can never really be sure that a given causal model represents
the entire cause-effect relations that hold among the given variables. 
In this context, we mention that the question of how to ensure that all hidden factors have been eliminated was one of the objections raised by the 
ancient C\={a}rv\={a}ka school of thought to inductive inference (see Page~15 of Chakrabarti (2010)). 

The modelling of causality using causal structures is a brilliant piece of deductive mathematical machinery. 
The {\em assumed} DAG and the probability distribution form the axioms of the deductive theory (along with the usual axioms of logic and probability). 
A theorem in the theory of causal structures is a statement about the DAG and the probability distribution and is proved using usual deductive mathematical
methods.  So, the remarkable theorem that it is possible to derive the effects of intervention without actual experiments is after all a deductive result. As such, it 
is non-ampliative in the sense that it does not convey more knowledge than what was already contained in the design of the causal structure. Statistical inference, on the 
other hand, is ampliative. So, at a fundamental level the applicability of modelling using causal structures to statistical inference is limited. 

As has been argued above, statistical inference is one form of inductive inference. In fact, as quoted above, Mahalanobis considers statistics to be {\em the} 
general method of inductive inference. As a method of inductive inference, broadly speaking one may consider one of the goals of statistical inference to derive
cause-effect relations. Statistical inference, however, avoids talking about cause-effect relations since statistical techniques cannot determine such relations. 
On the other hand, causal structures modelling starts out by postulating causal relations and the associated probability distributions. So, considering statistical inference 
to be an exercise in causal modelling is like postulating what one wishes to derive. In this context, it may be noted that Russell had raised objections
about ``introducing entities with implicit definitions, that is, as being those things that obey certain axioms or `postulates'.'' (Linsky (2019))

An important aspect of any inductive inference is non-monotonicity. In the context of SP, this issue has been discussed in Section~\ref{sec-data}. We note that
Imbens and Rubin (1995) remark that ``the monotonicity assumption is difficult to represent in a graphical model''.

\section{To Aggregate or Not to Aggregate? \label{sec-resolution} }

In Section~\ref{sec-condn}, we have provided several sufficient conditions for SP not to hold. These are conditions on the sub-populations. If any of these conditions hold,
then the paradox will not arise. Among these sufficient conditions are two homogeneity conditions, one of which is a slight modification of a condition due to Mittal (1991)
and the other is a new condition, namely the WORH condition. The sufficient conditions in general and the homogeneity conditions in particular, are not
necessary. Consequently, the question arises as to whether one should accept the conclusion of the aggregated table when the homogeneity conditions do not hold.
If one takes the position that the data are acceptable only if the homogeneity condition holds, then there will not be any paradox and effectively this
would make the resolution of the paradox very simple. However, if a homogeneity condition does not hold, then this does not necessarily mean that aggregation should 
not be done. There are scenarios, where the homogeneity conditions do not hold, yet, it is the aggregated table which is the correct answer. 
From Section~\ref{sec-ex}, one may note such examples arising in Simpson (1951) and the agricultural example in Lindley and Novick (1981). 
Nevertheless, the non-satisfaction of the homogeneity conditions does raise a red flag. In this context, we recall the caveat 
from Mittal (1991) who recommended that if the homogeneity condition she identified is not satisfied, then ``a closer look at the data is warranted before amalgamation ... 
Such a look may reveal some hidden factors in the data that might make amalgamation unwise.''

Simpson's paradox has been investigated from the empirical perspective by Spanos (2020) in the parametric framework using the notion of statistical
mis-specification. Homogeneity considerations arise in this investigation and Spanos (2020) points out that the untrustworthy associations in 
the Cohen and Nagel (1934) example is due to the fact that ``the two populations are not homogenous.'' Several other examples are also explained by showing that
aggregation leads to unreliable models. We note, on the other hand, that Spanos (2020) does not appear to provide an example, where homogeneity conditions are
violated and yet the aggregated data provides the correct answer. As mentioned above, such an example is the agricultural data example in 
Lindley and Novick (1981). It would be of interest to know whether the approach suggested by Spanos (2020) can explain this (and similar) examples.


From the discussion following Definition~\ref{defn-cont-tab}, we have an instance of the paradox if 
$\Pr[Y|X,M]>\Pr[Y|X,\overline{M}]$, $\Pr[Y|\overline{X},M]>\Pr[Y|\overline{X},\overline{M}]$, but $\Pr[Y|X]<\Pr[Y|\overline{X}]$. One way to understand how this can 
happen is to use the conditional probability expansions of $\Pr[Y|X]$ and $\Pr[Y|\overline{X}]$ as follows. 
\begin{eqnarray*}
	\Pr[Y|X] & = & \Pr[Y|X,M]\Pr[M|X] + \Pr[Y|X,\overline{M}]\Pr[\overline{M}|X], \\
	\Pr[Y|\overline{X}] & = & \Pr[Y|\overline{X},M]\Pr[M|\overline{X}] + \Pr[Y|\overline{X},\overline{M}]\Pr[\overline{M}|\overline{X}].
\end{eqnarray*}
So, both $\Pr[Y|X]$ and $\Pr[Y|\overline{X}]$ can be seen as weighted sums, where the weights for $\Pr[Y|X]$ are $\Pr[M|X]$ and $\Pr[\overline{M}|X]$ and the
weights for $\Pr[Y|\overline{X}]$ are $\Pr[M|\overline{X}]$ and $\Pr[\overline{M}|\overline{X}]$. Since the weights corresponding to 
$\Pr[Y|X,M]$ and $\Pr[Y|\overline{X},M]$ are not equal and similarly, the weights corresponding to $\Pr[Y|X,\overline{M}]$ and $\Pr[Y|\overline{X},\overline{M}]$
are not equal, the inequalities $\Pr[Y|X,M]>\Pr[Y|X,\overline{M}]$ and $\Pr[Y|\overline{X},M]>\Pr[Y|\overline{X},\overline{M}]$ are not preserved on aggregation.

The above provides a simple mathematical explanation of the paradox and has been pointed out by Hand (1994) who also mention that this is an ``inadequate
resolution.'' In the words of Hand (1994), ``the issue is whether or not ... to ask a conditional question and, in fact, the issue of whether or not to
condition is ubiquitous. ... Neither analysis is right and the other wrong -- it depends on what we want to find out.''
The issue of formulating separate questions whose answers are provided by the disaggregated and the aggregated tables is another dimension of the 
typical `What to do?' question that an investigator has to contend with when faced with an instance of SP.
We may take this view and apply it to some of the examples in Section~\ref{sec-ex}. 
\begin{enumerate}
\item Consider the example of Gardner (1976). If we wish to maximise the probability of drawing a blue chip, then if the chips are in the disaggregated format, we should 
draw from the black hats, while if the chips are in the aggregated format, then we should draw from the grey hat. Since, the chips will be in either the disaggregated or 
the aggregated format, we have a complete resolution of the paradox in this particular case.
\item Consider the agricultural example of Lindley and Novick (1981). If both short and tall varieties are available in adequate quantity, then it is better to
ignore this aspect and plant only the white variety. On the other hand, if only the short variety, or, only the tall variety is available, then clearly
one should go by the disaggregated table and plant the black variety. So, depending on the context, both the disaggregated and the aggregated tables provide
useful information. 
\end{enumerate}
The question arises as to whether the above kind of analysis can be applied to all examples of SP? It seems that there are difficulties in doing so. Consider
the Lindley and Novick (1981) example of treatment described in Section~\ref{sec-ex}. In this example, the moot question is whether to apply the treatment or not, to
which the answer is provided by the disaggregated tables; it seems difficult to come up with a question whose answer is given by the aggregated table. 

In Pearl (2014), it has been mentioned that ``in certain models the correct answer may not lie in either the disaggregated or the aggregated data.'' An 
explanation for such a situation has been provided in terms of the backdoor criterion. 
We would like to highlight that there is no known example of SP in the literature where neither the disaggregated nor the aggregated data provide the correct answer. 
This raises the question of whether this possibility is only a mathematical artifice, or, whether such cases can really arise in practice?
Further, we note that the possibility that for a particular example, both the disaggregated and the aggregated data could be correct depending on the context, does not 
seem to have been considered by Pearl.


\section{Broader Perspective: Logic, Probability, and Statistics \label{sec-broad} }

Issues involving causation in statistics date to more than one hundred year ago. There are at least four possible positions regarding the role and status of
causality in statistics. They are (i) radical causal skepticism, 
(ii) modest causal skepticism, (iii) causal/statistical compatibilism, and (iv) full-blown adherence to causality. Radical skepticism about causation is aptly 
captured by Karl Pearson for whom correlation is all there is. For him, there is no room for causation in statistics (See Aldrich, 1995 for a nuanced approach to Pearson). 
Causation skepticism dates to Pearson (1910) when Fisher's theory of randomization had not yet developed. 
In his words, (1910 quoted in Aldrich, 1995.)\footnote{In one sense, causality existed in statistics before 1950's. The role of randomization is where causality arises 
after Pearson left. It was precisely in the hand of Fisher that the theory was developed. There was an intellectual feud between Fisherians and Wright about causality and 
the role of path-diagrams in data analysis. It resulted in a sustained neglect of the appreciation of path analysis in the statistical community (see, Pearl and Mackenzie, 2018).}

\begin{quotation}
``It is this conception of correlation between two occurrences embracing all relationships from absolute independence to complete dependence, which is the wider 
	category by which we have to replace the old idea of causation (Pearson, 1910, p.157).''
\end{quotation}

No wonder that Russell, being an admirer of Pearson's work, shares his causal skepticism in science. Pearl calls Pearson ``causality's worst adversary.'' (Pearl, 2009, P.105.) 

A weaker position about statisticians' stance toward causal language can be attributed to Novick and Lindley's influential paper in which they address Simpson's Paradox 
from a Bayesian statistic. For them, informal causal talk is possible, but, according to them, causation is not a well-defined term to be presented in statistics. They write, 
``One possibility would be to use the language of causation .... We have not chosen to do this; nor to discuss causation, because the concept, although widely used, does not 
seem to be well-defined (Novick and Lindley, 1981).'' Their view on causation in statistics resembles Earman's oft-quoted expressions about ``the wooliness''
of causal notion in physics.
We call their position, modest causal skepticism. Causal theorists such as Glymour, Meek, and Pearl took them to task 
for trying hard to avoid causal language when clearly, according to these proponents of causal theories, causal notions are brought in their paper from the backdoor 
(because of the use of exchangeability). Elsewhere, Pearl has gone further than they and criticized the statistical community with a blistering language (Pearl, 2009, P.177.).

\begin{quotation}
``Simpson's paradox helps us to appreciate both the agony and the achievement of this generation of statisticians. Driven by healthy intuition, yet culturally forbidden from 
	admitting it and mathematically disabled from expressing it, they managed nevertheless to extract meaning from dry tables and to make statistical methods the standard 
	in the empirical sciences. But the spice of Simpson's paradox turned out to be nonstatistical (i.e., {\em causal}) after all.'' (Emphasis is ours.)
\end{quotation}

Those who provide a full-blown defense of causality in statistical side are Pearl on the one hand, and Spirtes, Glymour, and Schienes on the other. Pearl is a strong proponent 
of causality. Pearl is interested in the ``Why Question'' where a connection needs to be established between a known cause and a known effect for a better grasp of the mechanism 
at play. He discusses in particular the role causality plays in scientific inference. This discussion is in terms of what he calls ``the ladder of causation'' consisting of three 
levels. The first level is the world of association in which statistics plays a role, and there is no causation properly so called. The second level is where intervention plays 
a role. We intervene in a system to know its causal mechanism with the help of a causal model. The third and final level in the ladder is the world of counterfactuals in which 
there are no data to offer an insight into counterfactual reasoning. We ask, ``What would be the case for y if x had happened?''\footnote{David Lewis has contributed considerably 
to our understanding of counterfactuals (Lewis, 1973.)}

Pearl is interested in understanding the foundation of a theory of causality. However, he is keener on developing a mathematical language to furnish the theoretical 
underpinning of causality. He thinks that statistical apparatus fails to capture this kind of counterfactual reasoning as statistics, according to him, involves and pivots 
only around data. With the information of an able investigator, we can propose a causal model to understand the mechanism of some phenomenon where data could provide some 
clue as to their mechanism, but not necessarily the full story of the phenomenon at stake.

The CMU philosophers/statisticians share Pearl's enthusiasm for the significance of causality in contemporary methodology. There is at least one difference between Pearl 
and the CMU philosophers and statisticians, however. Pearl is not interested in the discovery of a cause, as he builds causal models by drawing on expert knowledge and 
available data. In contrast, CMU philosophers/statisticians are interested in causal discovery. To this end, they suggest a subject-matter-neutral automated causal inference 
engine that provides causal relationships among variables from observational data using information on their probabilistic correlations and assumptions about their causal structure.  
Both causal accounts are anti-reductionists in the sense that a causal relationship among events cannot be defined in terms of any other concepts in addition to the fact that 
data cannot provide a full narrative about the causal process. According to them, we are required to have causal assumptions to begin with to make causal inferences.

We reconstructed Pearl's strongest stance toward causality including his critical stance against statistical community for disregarding the need for 
casual language in providing a theory of scientific inference. However, he has changed his rhetoric and become less brutal about statistics being able to express causal language. 
Presumably, the most weakened stance he has taken toward statistics is to enrich statistics when it needs causal language as well as employing statistical language to provide a 
theory of causal inference. We call this position of the co-existence of statistics with causality, the causality-statistics compatibilism.

Pearl writes ``causation is not merely an aspect of statistics; it is an addition to statistics, an enrichment that allows statistics to uncover workings of the world that 
traditional methods cannot.'' (Pearl et al. 2016.) So, the analysis of what Pearl and his co-authors provide is that investigators need to know how and why causes influence their 
effects by analyzing data in a study. An analysis of the data helps appreciate why a cause that holds in one context may not hold in another. Presumably Pearl's view is that 
generalizability across domains is motivated by how science works, and not just from how statistics functions. The title of their book {\em Causal 
Inference in Statistics} is suggestive. It captures a mutual feedback process that causal inference can incorporate statistics, and the converse is the case; yet they are 
distinct disciplines (Pearce and Lawlor, 2016).

To summarize, there are at least four accounts of the role of statistics in causality debate. They are (i) radical causal skepticism, (ii) modest causal skepticism, 
(iii) causal/statistical compatibilism, and finally, the full-blown causal theory of scientific inference.

Keeping in mind our interest in Simpson's Paradox, we will only mention 
two extreme stances toward causality. They are the Russell-Pearson position where causality is completely eschewed, and the Pearl-CMU view in which causality is embraced 
whole-heartedly. Our Simpson's Paradox-centered stance concerning causality does not adopt either of the extreme positions. 

The underlying theme of this paper is to contend that we need causal language combined with statistical/probabilistic tools in which philosophy helps sort out among different 
types of questions. We raised at least three questions regarding the paradox: (i) Why is Simpson’s Paradox paradoxical? (ii) What are the conditions for generating the paradox, 
and (iii) What-to-do when confronted with a Simpson's Paradoxical type situation? We disagreed with the causal resolution regarding the surprising element in the paradox. While 
discussing different definitions of Simpson's Paradox, we showed that logic/probability coupled with statistical language devoid of causal intuition is adequately powerful to 
represent the paradox. Therefore, Simpson's Paradox has nothing to do with causal considerations in addressing the first two questions. 
For the what-to-do question, causal considerations, along with other statistical tools, are important in making inductive inferences.
Therefore, we are, in a sense, compatibilists regarding the co-existence of statistics and causality to be able to provide a sound theory of scientific inference.

\section*{Acknowledgement} We are grateful to Gordon Brittan for providing helpful comments on an earlier draft.

\section*{References}

J. von Kugelgen, L. Gresele and B. Scholkopf, "Simpson's Paradox in COVID-19 Case Fatality Rates: A Mediation Analysis of Age-Related Causal Effects" in IEEE Transactions on Artificial Intelligence, vol. 2, no. 01, pp. 18-27, 2021.
doi: 10.1109/TAI.2021.3073088
\begin{description}
	\item{Armistead, T.~W.:} Resurrecting the third variable: a critique of Pearl's causal analysis of Simpson's paradox. {\em The American Statistician}, 68(1), 1--7 (2014).
	\item{Aldrich, A.:} Correlations genuine and spurious in Pearson and Yule. {\em Statistical Science}, 10(4), 364--376 (1995).
	\item{Bandyopadhyay, P.~S., Nelson, D., Greenwood, M., Brittan, G., and Berwald, J.:} The logic of Simpson's paradox. {\em Synthese}, 181, 185--208 (2011). 
	\item{Bickel, P.~J., Hammel, E.~A., and O'Connell, J.~W.:} Sex bias in graduate admissions. {\em Science}, 187, 398--404 (1975).
	\item{Blyth, C.:} On Simpson's paradox and the sure-thing principle. {\em Journal of the American Statistical Association}, 67, 364--366 (1972).
	\item{Chakrabarti, K.~K.:} {\em Classical Indian philosophy of induction: the Nyaya viewpoint}. Lexington Books, 2010.
	\item{Chakravarty, S. and Sarkar, P.:} An Inequality Paradox: Relative versus Absolute Indices? {\em Metron}, 79, 241--254 (2021).
	\item{Cohen, M. and Nagel, E.:} {\em An Introduction to Logic and the Scientific Method}. Harcourt, Brace and Company, New York, 1934.
	\item{Cox, D.R., and Wermuth, N.:} Discussion of `Causal diagrams for empirical research' by J.~Pearl. {\em Biometrika}, 82(4), 688--689 (1995).
	\item{Daudt, H., and D.W. Rae, D.W.:} The Ostrogorski paradox: a peculiarity of compound majority decision, {\em European Journal of Political Research}, 4(4), 391-399 (1976).
	\item{Fisher, R.:} Statistical methods and scientific induction. {\em Journal of the Royal Statistical Society, Series B (Methodological)}, 17(1), 69--78 (1955).
	\item{Fitelson, B.:} (2017). Confirmation, causation, and Simpson's paradox. Episteme, 14(3), 297--309 (2017). 
	\item{Freeman, D.:} Discussion of `Causal diagrams for empirical research' by J.~Pearl. {\em Biometrika}, 82(4), 692--693 (1995).
	\item{Gardner, M.:} ``Mathematical Games.'' {\em Scientific American}, 234(3), 119--124 (1976).
	\item{Good,~I.~J., and Mittal.~Y.:} The amalgamation and geometry of two-by-two contingency tables. {\em The Annals of Statistics}, 15(7), 694--711 (1987).
	\item{Hand,~D.~J.:} Deconstructing statistical questions. {\em Journal of the Royal Statistical Society, Series A}, 157(3), 317--356 (1994).
	\item{Imbens,~G.~W., and Rubin,~D.B.:} Discussion of `Causal diagrams for empirical research' by J.~Pearl. {\em Biometrika}, 82(4), 694--695 (1995).
	\item{Kornhauser, L. A., and Sager, L.G.:} (1986). Unpacking the court, {\em Yale Law Journal}, 96(1), 82--117, (1986).
	\item{Kugelgen, J. von, Gresele, L. and Scholkopf, B.:} Simpson's paradox in COVID-19 case fatality rates: a mediation analysis of age-related causal effects. 
		{\em IEEE Transactions on Artificial Intelligence}, 2(1), 18--27, 2021.  
	\item{Levi, I.:} Inductive inference as ampliative and non monotonic reasoning. {\em TARK'05: Proceedings of the 10th conference on Theoretical Aspects of Rationality 
		and Knowledge}, 177--192 (2005).
	\item{Lewis, D.:} Causation. {\em Journal of Philosophy}, 70, 556--567 (1973).
	\item{Lewis, D.:} Philosophical Papers. Vol. II. Oxford University Press. 1986.
	\item{Lindley, D., and Novick, M.:} The role of exchangeability in inference. {\em The Annals of Statistics}, 9, 45--58 (1981).
	\item{Linsky, B.:} {\em Logical Constructions}. The Stanford Encyclopedia of Philosophy, \url{https://plato.stanford.edu/entries/logical-construction/} (2019).
	\item{List, C., and Pettit, P.:}. Aggregating sets of judgments: An impossibility result. {\em Economics and Philosophy}, 18(1), 89--110, (2002).
	\item{Mahalanobis, P.C.:} Why statistics? {\em Sankhy\={a}: The Indian Journal of Statistics}, 10(3), 195--228, (1950).
	\item{Mayo,~D.G., and Cox,~D.~R.:} Frequentist statistics as a theory of inductive inference, {\em Institute of Mathematical Statistics Lecture Notes - Monograph Series},
		49, 77--97, (2006).
	\item{Mittal, Y.:} Homogeneity of subpopulations and Simpson's paradox. {\em Journal of the American Statistical Association}, 86, 167--172 (1991). 
	\item{Pagano, M., and Gauvreau, K.:} {\em Principles of Biostatistics}, second edition. Duxbury, 1993.
	\item{Pearce, N, and Lawlor, D.:} Causal inference -- so much more than statistics. {\em International Journal of Epidemiology}, 1895--1903, (2016).
	\item{Pearl, J. and Mackenzie. D.:} {\em The book of why}. Basic Books, New York (2018).
	\item{Pearl J, Glymour M, and Jewell N.:} Causal inference in statistics: a primer. {\em Wiley}, UK, (2016).
	\item{Pearl, J.:} Comment: understanding Simpson's paradox. {\em The American Statistician}, 68(1), 8--13 (2014).
	\item{Pearl, J.:} {\em Causality: models, reasoning and inference}, second edition. Cambridge University Press, 2009.
	\item{Pearl, J.:} Causal diagrams for empirical research. {\em Biometrika}, 82(4), 669--688 (1995).
	\item{Pearson, K., Lee, A., and Bramley-Moore, L.:} Genetic (reproductive) selection: inheritance of fertility in man, and of fecundity in thoroughbred
		race horses. {\em Philosophical Transactions of the Royal Society of London, Series A}, 192, 257--330 (1899).
	\item{Quine, W. V.:} The ways of paradox and other essays. Harvard University, Press. 1976.
	\item{Rao,~C.~R.:} {\em Statistics and truth: putting chance to work}, second edition, World Scientific, 1997.
	\item{Russell, B.:} On the notion of cause. In the Collected Papers of Bertrand Russell, vol 12, London: Routledge, 1992, 193--210. Originally published in 1913.
	\item{Simpson, E.:} The interpretation of interaction in contingency tables. {\em Journal of the Royal Statistical Society, Series B}, 13, 238--241 (1951).
	\item{Spanos, A.:} Yule-Simpson's paradox: the probabilistic versus the empirical conundrum. {\em Statistical Methods \& Applications}, 
		\url{https://doi.org/10.1007/s10260-020-00536-4} (2020).
	\item{Selvitella, A.:} The ubiquity of the Simpson's paradox. Journal of Statistical Distributions and Applications 4, 1--16 (2017).
	\item{Selvitella, A.:} The Simpson's paradox in quantum mechanics. Journal of Mathematical Physics, 58 (2017a).
	\item{Sober, E.:} Core Questions in Philosophy: A Text Readings. 8th Edition. 2021. What is Philosophy?
	\item{Spirties, P., Glymour, C., and Scheines, R.:} {\em Causation, prediction, and search}. MIT Press, Cambridge (2000).
	\item{Sprenger, J. and Weinberger, N.:} Simpson's Paradox, The Stanford Encyclopedia of Philosophy (Summer 2021 Edition), Edward N. Zalta (ed.), 
		\url{https://plato.stanford.edu/archives/sum2021/entries/paradox-simpson/}.
	\item{Yule, G.:} Notes on the theory of association of attributes in statistics. {\em Biometrika}, 2, 121--134 (1903).
\end{description}

\appendix

\section{Proofs \label{sec-proof} }

\subsection{Proof of Theorem~\ref{thm-possible} \label{subsec-pf-thm-possible} }

\begin{proof}
The result is proved by showing an example for each of the cases. Note that for $1\leq i\leq 13$, Case number $(28-i)$
transforms to Case number $i$ by interchanging the rows of $T_1$, $T_2$ and $T_1+T_2$. 
Further, by interchanging the tables $T_1$ and $T_2$, Cases~4,~5 and~6 get transformed to Cases~10,~11 and~12 respectively. 
In view of these observations, it is sufficient to provide examples for Cases~1 to~9 and~13,~14. Figure~\ref{fig-ex} provides these examples. 

\begin{figure}
	\caption{Examples of $T_1$, $T_2$ and $T_1+T_2$ satisfying Cases~1 to~14 of Figure~\ref{fig-reln}. \label{fig-ex} }
\begin{center}
\begin{tabular}{|l|c|c|c|}
\hline
	Case & $T_1$ & $T_2$ & $T_1+T_2$ \\ \hline\hline
	1 & $\begin{array}{|c|c|} \hline 3 & 1 \\ \hline 1 & 1 \\ \hline \end{array}$ 
	& $\begin{array}{|c|c|} \hline 3 & 1 \\ \hline 1 & 1 \\ \hline \end{array}$ 
	& $\begin{array}{|c|c|} \hline 6 & 2 \\ \hline 2 & 2 \\ \hline \end{array}$ \\ \hline
	2 & $\begin{array}{|c|c|} \hline 3 & 1 \\ \hline 6 & 3 \\ \hline \end{array}$ 
	& $\begin{array}{|c|c|} \hline 5 & 7 \\ \hline 1 & 4 \\ \hline \end{array}$ 
	& $\begin{array}{|c|c|} \hline 8 & 8 \\ \hline 7 & 7 \\ \hline \end{array}$ \\ \hline
	3 & $\begin{array}{|c|c|} \hline 5 & 3 \\ \hline 10 & 10 \\ \hline \end{array}$ 
	& $\begin{array}{|c|c|} \hline 1 & 19 \\ \hline 1 & 20 \\ \hline \end{array}$ 
	& $\begin{array}{|c|c|} \hline 6 & 22 \\ \hline 11 & 30 \\ \hline \end{array}$ \\ \hline
	4 & $\begin{array}{|c|c|} \hline 3 & 1 \\ \hline 1 & 1 \\ \hline \end{array}$ 
	& $\begin{array}{|c|c|} \hline 2 & 2 \\ \hline 3 & 3 \\ \hline \end{array}$ 
	& $\begin{array}{|c|c|} \hline 5 & 3 \\ \hline 4 & 4 \\ \hline \end{array}$ \\ \hline
	5 & $\begin{array}{|c|c|} \hline 2 & 8 \\ \hline 1 & 5 \\ \hline \end{array}$ 
	& $\begin{array}{|c|c|} \hline 4 & 4 \\ \hline 3 & 3 \\ \hline \end{array}$ 
	& $\begin{array}{|c|c|} \hline 6 & 12 \\ \hline 4 & 8 \\ \hline \end{array}$ \\ \hline
	6 & $\begin{array}{|c|c|} \hline 2 & 1 \\ \hline 3 & 2 \\ \hline \end{array}$ 
	& $\begin{array}{|c|c|} \hline 3 & 3 \\ \hline 1 & 1 \\ \hline \end{array}$ 
	& $\begin{array}{|c|c|} \hline 5 & 4 \\ \hline 4 & 3 \\ \hline \end{array}$ \\ \hline
	7 & $\begin{array}{|c|c|} \hline 10 & 5 \\ \hline 3 & 5 \\ \hline \end{array}$ 
	& $\begin{array}{|c|c|} \hline 1 & 2 \\ \hline 1 & 1 \\ \hline \end{array}$ 
	& $\begin{array}{|c|c|} \hline 11 & 7 \\ \hline 4 & 6 \\ \hline \end{array}$ \\ \hline
	8 & $\begin{array}{|c|c|} \hline 2 & 1 \\ \hline 1 & 3 \\ \hline \end{array}$ 
	& $\begin{array}{|c|c|} \hline 1 & 5 \\ \hline 1 & 1 \\ \hline \end{array}$ 
	& $\begin{array}{|c|c|} \hline 3 & 6 \\ \hline 2 & 4 \\ \hline \end{array}$ \\ \hline
	9 & $\begin{array}{|c|c|} \hline 3 & 2 \\ \hline 1 & 1 \\ \hline \end{array}$ 
	& $\begin{array}{|c|c|} \hline 1 & 5 \\ \hline 3 & 2 \\ \hline \end{array}$ 
	& $\begin{array}{|c|c|} \hline 4 & 7 \\ \hline 4 & 3 \\ \hline \end{array}$ \\ \hline
	13 & $\begin{array}{|c|c|} \hline 2 & 2 \\ \hline 3 & 3 \\ \hline \end{array}$ 
	& $\begin{array}{|c|c|} \hline 1 & 2 \\ \hline 2 & 4 \\ \hline \end{array}$ 
	& $\begin{array}{|c|c|} \hline 3 & 4 \\ \hline 5 & 7 \\ \hline \end{array}$ \\ \hline
	14 & $\begin{array}{|c|c|} \hline 2 & 2 \\ \hline 3 & 3 \\ \hline \end{array}$ 
	& $\begin{array}{|c|c|} \hline 2 & 2 \\ \hline 2 & 2 \\ \hline \end{array}$ 
	& $\begin{array}{|c|c|} \hline 4 & 4 \\ \hline 5 & 5 \\ \hline \end{array}$ \\ \hline
\end{tabular}
\end{center}
\end{figure}
\end{proof}

\subsection{Proof of Theorem~\ref{thm-SP-nec-cond} \label{subsec-pf-thm-SP-nec-cond} }
\begin{proof} 
	Consider the first point. 
	\begin{itemize}	
	\item We show that $A_1>C_1$, $A_2>C_2$ and $A_1=A_2$ together imply $\mu>\nu$ and so $\sym{SP}_1$ cannot hold. 
	Let $k=A_1=A_2$. Then $a_1=k\alpha_1$ and $a_2=k\alpha_2$ and $\mu = (a_1+a_2)/(\alpha_1+\alpha_2)=k$.
	From $k=A_1>C_1$, we have $k\gamma_1>c_1$ and from $k=A_2>C_2$, we have $k\gamma_2>c_2$. So, $c_1+c_2<k(\gamma_1+\gamma_2)$ which shows that
	$\mu=k>(c_1+c_2)/(\gamma_1+\gamma_2)=\nu$.
	\item A similar argument shows that $A_1<C_1$, $A_2<C_2$ and $A_1=A_2$ together imply $\mu<\nu$ and so $\sym{SP}_2$ cannot hold.
	\end{itemize}
	So, if $A_1=A_2$, then neither $\sym{SP}_1$ nor $\sym{SP}_2$ hold and so SP does not hold. Contrapositively, if SP holds, then $A_1\neq A_2$.

Consider the second point.
	\begin{itemize}
	\item We show that $A_1>C_1$, $A_2>C_2$ and $C_1=C_2$ together imply $\mu>\nu$ and so $\sym{SP}_1$ cannot hold. 
	Let $s=C_1=C_2$. Then $c_1=s\gamma_1$ and $c_2=s\gamma_2$ and $\nu = (c_1+c_2)/(\gamma_1+\gamma_2)=s$.
	From $A_1>C_1=s$, we have $a_1>s\alpha_1$ and from $A_2>C_2=s$, we have $a_2>s\alpha_2$. So, $a_1+a_2>s(\alpha_1+\alpha_2)$ which shows that
	$\mu=(a_1+a_2)/(\alpha_1+\alpha_2)>s=\nu$.
	\item A similar argument shows that $A_1<C_1$, $A_2<C_2$ and $A_1=A_2$ together imply $\mu<\nu$ and so $\sym{SP}_2$ cannot hold.
	\end{itemize}
	So, if $C_1=C_2$, then neither $\sym{SP}_1$ nor $\sym{SP}_2$ hold and so SP does not hold. Contrapositively, if SP holds, then $C_1\neq C_2$.

Consider the third point.
        \begin{itemize}
	\item We show that $A_1>C_1$, $A_2>C_2$, $\alpha_1=\gamma_1$ and $\alpha_2=\gamma_2$ together imply $\mu>\nu$ and so $\sym{SP}_1$ cannot hold. 
	From $A_1>C_1$ and $\alpha_1=\gamma_1$, we have $a_1>c_1$; from $A_2>C_2$ and $\alpha_2=\gamma_2$, we have $a_2>c_2$. So, $a_1+a_2>c_1+c_2$. From
	$\alpha_1=\gamma_1$ and $\alpha_2=\gamma_2$, we have $\alpha_1+\alpha_2=\gamma_1+\gamma_2=x$ (say). Then 
	$\mu=(a_1+a_2)/(\alpha_1+\alpha_2)=(a_1+a_2)/x>(c_1+c_2)/x=(c_1+c_2)/(\gamma_1+\gamma_2)=\nu$.
	\item A similar argument shows that $A_1<C_1$, $A_2<C_2$, $\alpha_1=\gamma_1$ and $\alpha_2=\gamma_2$ together imply $\mu<\nu$ and so $\sym{SP}_2$ cannot hold.
	\end{itemize}
	So, if $\alpha_1=\gamma_1$ and $\alpha_2=\gamma_2$, then neither $\sym{SP}_1$ nor $\sym{SP}_2$ hold and so SP does not hold. Contrapositively, if SP holds, 
	then $\alpha_1\neq \gamma_1$ or $\alpha_2\neq \gamma_2$.
\end{proof}

\subsection{Proof of Theorem~\ref{thm-SP-basic} \label{subsec-pf-thm-SP-basic} }
\begin{proof} 
For $i=1,2$, note that $\kappa(T_i)>1$, i.e., $a_id_i>b_ic_i$ holds if and only if 
\begin{eqnarray}\label{eqn-tmp1}
\frac{a_i}{a_i+b_i} & > & \frac{c_i}{c_i+d_i}.
\end{eqnarray}
Similarly, $\kappa(T_1+T_2)>1$, i.e., $(a_1+a_2)(d_1+d_2) > (b_1+b_2)(c_1+c_2)$ holds if and only if 
\begin{eqnarray}\label{eqn-tmp2}
\frac{a_1+a_2}{a_1+a_2+b_1+b_2} & > & \frac{c_1+c_2}{c_1+c_2+d_1+d_2}.
\end{eqnarray}
So, $\sym{SP}_1$ is equivalent to the following: $ \left(\kappa(T_1) > 1\right) \wedge \left(\kappa(T_2) > 1\right) \wedge \neg\left(\kappa(T_1+T_2) > 1\right).  $
Similarly, $\sym{SP}_2$ is equivalent to the following: $ \left(\kappa(T_1) < 1\right) \wedge \left(\kappa(T_2) < 1\right) \wedge \neg\left(\kappa(T_1+T_2) < 1\right). $
Consequently, we have that $\sym{SP}_1\vee \sym{SP}_2$ is equivalent to~\eqref{eqn-SP-OR}.
\end{proof}

\subsection{Proof of Theorem~\ref{thm-SP-char} \label{subsec-pf-thm-SP-char} }
\begin{proof} From~(\ref{eqn-SP-OR}), Simpson's paradox does not hold if and only if the following condition holds.
$$
\neg\left(\left(\kappa(T_1) > 1\right) \wedge \left(\kappa(T_2) > 1\right) \wedge \neg\left(\kappa(T_1+T_2) > 1\right)\right) \\
	\mbox{ and }
\neg\left(\left(\kappa(T_1) < 1\right) \wedge \left(\kappa(T_2) < 1\right) \wedge \neg\left(\kappa(T_1+T_2) < 1\right)\right).
$$
Using basic logical equivalences, we obtain
\begin{eqnarray*}
\lefteqn{\neg\left(\left(\kappa(T_1) > 1\right) \wedge \left(\kappa(T_2) > 1\right) \wedge \neg\left(\kappa(T_1+T_2) > 1\right)\right)} \\
& \equiv & 
\neg\left(\kappa(T_1) > 1\right) \vee \neg\left(\kappa(T_2) > 1\right) \vee \left(\kappa(T_1+T_2) > 1\right) \\
& \equiv & 
\neg\left(\left(\kappa(T_1) > 1\right) \wedge \left(\kappa(T_2) > 1\right)\right) \vee \left(\kappa(T_1+T_2) > 1\right) \\
& \equiv &
\left(\left(\kappa(T_1) > 1\right) \wedge \left(\kappa(T_2) > 1\right)\right) \mbox{ implies } \left(\kappa(T_1+T_2) > 1\right).
\end{eqnarray*}
Similarly, one can argue that the condition
	$\neg\left( \left(\kappa(T_1) < 1\right) \wedge \left(\kappa(T_2) < 1\right) \wedge \neg\left(\kappa(T_1+T_2) < 1\right)  \right)$ is equivalent
	to 
	$\left(\kappa(T_1) < 1\right) \wedge \left(\kappa(T_2) < 1\right) \mbox{ implies } \kappa(T_1+T_2) < 1$.
\end{proof}


\subsection{Proof of Theorem~\ref{thm-suff-m91} \label{subsec-pf-thm-suff-m91} }
\begin{proof}
	We need to show that if any one of the conditions given in Definition~\ref{defn-homogeneity} hold, then SP does not occur. We show this for the first condition. 
	The proofs for the other conditions are similar. 

	Consider the first condition, namely, $\max(a_1/b_1,a_2/b_2)<\min(c_1/d_1,c_2/d_2)$. This is equivalent to the following four conditions: 
	$a_1d_1<b_1c_1$, $a_2d_2<b_2c_2$, $a_1d_2<c_2b_1$ and $a_2d_1<c_1b_2$.
	From $a_1d_1<b_1c_1$ and $a_2d_2<b_2c_2$ we have $\kappa(T_1)<1$ and $\kappa(T_2)<1$. So, the antecedent of~\eqref{eqn-notSP1} is false and so~\eqref{eqn-notSP1}
	is true. The antecedent of~\eqref{eqn-notSP2}, on the other hand, is true. The consequent of~\eqref{eqn-notSP2} is $\kappa(T_1+T_2)<1$. 
	The last condition is equivalent to $(a_1+a_2)(d_1+d_2) < (b_1+b_2)(c_1+c_2)$. From the four inequalities indentified above, it follows that 
	each of the cross product term in $(a_1+a_2)(d_1+d_2)$ is less than a corresponding and unique cross-product term in $(b_1+b_2)(c_1+c_2)$. So, 
	we have $\kappa(T_1+T_2)<1$ and therefore~\eqref{eqn-notSP2} is true. Consequently, SP does not hold. 
\end{proof}

\subsection{Proof of Theorem~\ref{thm-WORH-notSP} \label{subsec-pf-thm-WORH-notSP} }
\begin{proof}
The WORH states that either $\kappa(T_1)=\kappa(T_1+T_2)$ or $\kappa(T_2)=\kappa(T_1+T_2)$. 
Suppose that $\kappa(T_1)=\kappa(T_1+T_2)$, the other case being similar. 

If either $\kappa(T_1)=1$ or $\kappa(T_2)=1$, then the antecedents of both~\eqref{eqn-notSP1} and~\eqref{eqn-notSP2} are false and consequently,
	both~\eqref{eqn-notSP1} and~\eqref{eqn-notSP2} are true. So, suppose that $\kappa(T_1)\neq 1$ and $\kappa(T_2)\neq 1$. 
	This leads to four cases, namely,
	Case-1: ($\kappa(T_1)>1$ and $\kappa(T_2)>1$); 
	Case-2: ($\kappa(T_1)>1$ and $\kappa(T_2)<1$); 
	Case-3: ($\kappa(T_1)<1$ and $\kappa(T_2)>1$); 
	Case-4: ($\kappa(T_1)<1$ and $\kappa(T_2)<1$).
	For Cases-2 and 3, the antecedents of both~\eqref{eqn-notSP1} and~\eqref{eqn-notSP2} are false and consequently,
        both~\eqref{eqn-notSP1} and~\eqref{eqn-notSP2} are true.
	Now Consider Case-1. The antecedent of~\eqref{eqn-notSP2} is false and so~\eqref{eqn-notSP2} is true. On the other hand, the antecedent of~\eqref{eqn-notSP1}
	is true. Using $\kappa(T_1)=\kappa(T_1+T_2)$, we have $\kappa(T_1+T_2)>1$ so that the consequent of~\eqref{eqn-notSP1} is also true. So,~\eqref{eqn-notSP1}
	is true. So, in Case-1, both~\eqref{eqn-notSP1} and~\eqref{eqn-notSP2} are true and therefore SP does not hold. The argument for Case-4 is similar.
\end{proof}

\subsection{Proof of Theorem~\ref{thm-odd-ratio-pos-assoc} \label{subsec-pf-thm-odd-ratio-pos-assoc} }

\begin{proof}
	The random variables $X$ and $Y$ are positively associated if $\Pr[Y|X]>\Pr[Y|\overline{X}]$, i.e., if $a/(a+b) > c/(c+d)$. The last condition
	holds if and only if $ad \geq bc$, i.e., if and only if $\kappa(T)>1$. The argument for the positive association of $X$ and $\overline{Y}$ is similar.
\end{proof}

\subsection{Proof of Theorem~\ref{thm-nec-cond} \label{subsec-pf-thm-nec-cond} }

\begin{proof}
To start with we translate the given conditions in terms of $a_i,b_i,c_i$ and $d_i$, $i=1,2$. 
	The conditions $(X\sim Y)|M$ and $(X\sim Y)|\overline{M}$ both hold, so using Theorem~\ref{thm-odd-ratio-pos-assoc}, we have
	$\kappa(T_i)>1$, $i=1,2$. It is given that SP holds for $(T_1,T_2)$. From Theorem~\ref{thm-SP-basic}, it follows that $\kappa(T_1+T_2)\leq 1$.
	Note that $\kappa(T_i)>1$ is equivalent to the condition $a_id_i>b_ic_i$ and $\kappa(T_1+T_2)\leq 1$ is equivalent to the condition
	$(a_1+a_2)(d_1+d_2)\leq (b_1+b_2)(c_1+c_2)$. It is possible to write these conditions in two different ways. 
	\ \\

	\noindent{\bf Case} $a_i/c_i>b_i/d_i$, $i=1,2$ and $(a_1+a_2)/(c_1+c_2) \leq (b_1+b_2)/(d_1+d_2)$: One may write 
	$(a_1+a_2)/(c_1+c_2) = \lambda (a_1/c_1) + (1-\lambda)(a_2/c_2)$ where $\lambda =c_1/(c_1+c_2)$ and
	$(b_1+b_2)/(d_1+d_2) = \mu (b_1/d_1) + (1-\mu)(b_2/d_2)$ where $\mu=d_1/(d_1+d_2)$. So, in particular $(a_1+a_2)/(c_1+c_2)$ lies in the
	interval defined by $a_1/c_1$ and $a_2/c_2$ and $(b_1+b_2)/(d_1+d_2)$ lies in the interval defined by $b_1/d_1$ and $b_2/d_2$.

	Since $a_i/c_i>b_i/d_i$, $i=1,2$, and $(a_1+a_2)/(c_1+c_2) \leq (b_1+b_2)/(d_1+d_2)$, we have 
	$\mu (b_1/d_1) + (1-\mu)(b_2/d_2) = (b_1+b_2)/(d_1+d_2) \geq (a_1+a_2)/(c_1+c_2) = \lambda (a_1/c_1) + (1-\lambda)(a_2/c_2)>\lambda (b_1/d_1) + (1-\lambda)(b_2/d_2)$. 
	So, $(\mu-\lambda)(b_1/d_1-b_2/d_2)>0$ and consequently, $\mu>\lambda$ or $\mu<\lambda$ according as $b_1/d_1>b_2/d_2$ or $b_1/d_1<b_2/d_2$; in particular, 
	neither $\mu=\lambda$ nor $b_1/d_1=b_2/d_2$ are possible. Using the definitions of $\mu$ and $\lambda$, we have 
	$c_1/c_2<d_1/d_2$ or $c_1/c_2>d_1/d_2$ according as $b_1/d_1>b_2/d_2$ or $b_1/d_1<b_2/d_2$.
	
	The condition $(a_1+a_2)/(c_1+c_2) \leq (b_1+b_2)/(d_1+d_2)$ forces the intervals $I=[b_1/d_1, a_1/c_1]$ and $J=[b_2/d_2, a_2/c_2]$ to be disjoint. To see
	this, first suppose that $a_2/c_2\geq a_1/c_1$. If the intervals $I$ and $J$ are not disjoint, then we have the relation
	$b_1/d_1 \leq b_2/d_2 \leq a_1/c_1 \leq a_2/c_2$. Since $(a_1+a_2)/(c_1+c_2)$ lies in $[a_1/c_1, a_2/c_2]$ and $(b_1+b_2)/(d_1+d_2)$ lies in
	$[b_1/d_1,b_2/d_2]$, it follows that $(a_1+a_2)/(c_1+c_2) \geq (b_1+b_2)/(d_1+d_2)$. Equality occurs if and only if 
	$(b_1+b_2)/(d_1+d_2)=b_2/d_2=a_1/c_1=(a_1+a_2)/(c_1+c_2)$, i.e., if and only if $b_1/d_1=b_2/d_2=a_1/c_1=a_2/c_2$. This in particular means that
	$b_1/d_1=a_1/c_1$. Since it is given that $a_1/c_1 > b_1/d_1$, equality cannot occur and so, $(a_1+a_2)/(c_1+c_2) > (b_1+b_2)/(d_1+d_2)$. The last condition
	contradicts $(a_1+a_2)/(c_1+c_2) \leq (b_1+b_2)/(d_1+d_2)$. Similarly, if $a_1/c_1\geq a_2/c_2$, then also we obtain a contradiction. This establishes
	that $I$ and $J$ are disjoint intervals. So, there are two possibilities: 
	(A) $b_1/d_1 < a_1/c_1 < b_2/d_2 < a_2/c_2$, or (B) $b_2/d_2 < a_2/c_2 < b_1/d_1 < a_1/c_1$. From the previous discussion, we obtain that 
	in (A), $c_1/c_2>d_1/d_2$ holds, while in (B) $c_1/c_2<d_1/d_2$ holds.
	\ \\

	\noindent{\bf Case} $a_i/b_i>c_i/d_i$, $i=1,2$ and $(a_1+a_2)/(b_1+b_2) \leq (c_1+c_2)/(d_1+d_2)$: As argued above, note that $(a_1+a_2)/(b_1+b_2)$ lies in the
	interval defined by $a_1/b_1$ and $a_2/b_2$ and similarly, $(c_1+c_2)/(d_1+d_2)$ lies in the interval defined by $c_1/d_1$ and $c_2/d_2$. Also, arguing as
	above, the condition $(a_1+a_2)/(b_1+b_2) \leq (c_1+c_2)/(d_1+d_2)$ forces the intervals $[c_1/d_1,a_1/b_1]$ and $[c_2/d_2,a_2/b_2]$ to be disjoint. This leads
	to the following two possibilities: (C) $c_1/d_1<a_1/b_1<c_2/d_2<a_2/b_2$, or (D) $c_2/d_2<a_2/b_2<c_1/d_1<a_1/b_1$. In (C), we have $c_1/c_2<d_1/d_2$, while
	in (D) we have $c_1/c_2>d_1/d_2$. 

	The combined effect of the above two cases is that if $\kappa(T_i)>1$, $i=1,2$ and $\kappa(T_1+T_2)\leq 1$, then either (A) and (C) both occur, or (B) and (D) both occur. 
	Since $(a_2+b_2)/(c_2+d_2)$ lies in $[b_2/d_2,a_2/c_2]$ and $(a_1+b_1)/(c_1+d_1)$ lies in $[b_1/d_1,a_1/c_1]$, 
	(B) implies $(a_2+b_2)/(c_2+d_2) < (a_1+b_1)/(c_1+d_1)$ which is equivalent to $M\sim X$. A similar reasoning shows that (D) implies $M\sim Y$. So, if (B) and (D) hold,
	then both $M\sim X$ and $M\sim Y$ hold. On the other hand, an analogous reasoning shows that if (A) and (C) hold, then both $M\sim \overline{X}$ and $M\sim \overline{Y}$
	hold.
\end{proof}

\subsection{Proof of Theorem~\ref{thm-togglings} \label{subsec-pf-thm-togglings} }
\begin{proof}
For a non-zero $\delta$ and a $2\times 2$ table $T$, by $\delta T$, we denote the $2\times 2$ table obtained from $T$ by multiplying each entry with $\delta$. 
It is easy to see that for $2\times 2$ tables $T_1$ and $T_2$, if $(T_1,T_2)$ is an instance of positive SP, then so is $(\delta T_1,\delta T_2)$. Similarly,
for positive alignment, negative SP and negative alignment.

Define
\begin{center}
\begin{tabular}{ccc}
\begin{minipage}{100pt}
\begin{eqnarray*}
	T_1^{(0)} & = &
\begin{array}{|c|c|}
\hline
5 & 3 \\ \hline
10 & 10 \\ \hline
\end{array}\ ,
\end{eqnarray*}
\end{minipage}
&
\begin{minipage}{100pt}
\begin{eqnarray*}
	T_2^{(0)} & = &
\begin{array}{|c|c|}
\hline
1 & 19 \\ \hline
1 & 20 \\ \hline
\end{array}\ ,
\end{eqnarray*}
\end{minipage}
&
\begin{minipage}{100pt}
\begin{eqnarray*}
	T_1^{(0)}+T_2^{(0)} & = &
\begin{array}{|c|c|}
\hline
6 & 22 \\ \hline
11 & 30 \\ \hline
\end{array}\ .
\end{eqnarray*}
\end{minipage} \\

\begin{minipage}{100pt}
\begin{eqnarray*}
	T_1^{(1)} & = &
\begin{array}{|c|c|}
\hline
10 & 6 \\ \hline
10 & 10 \\ \hline
\end{array}\ ,
\end{eqnarray*}
\end{minipage}
&
\begin{minipage}{100pt}
\begin{eqnarray*}
	T_2^{(1)} & = &
\begin{array}{|c|c|}
\hline
1 & 19 \\ \hline
1 & 20 \\ \hline
\end{array}\ ,
\end{eqnarray*}
\end{minipage}
&
\begin{minipage}{100pt}
\begin{eqnarray*}
	T_1^{(1)}+T_2^{(1)} & = &
\begin{array}{|c|c|}
\hline
11 & 25 \\ \hline
11 & 30 \\ \hline
\end{array}\ .
\end{eqnarray*}
\end{minipage} \\

\begin{minipage}{100pt}
\begin{eqnarray*}
	T_1^{(2)} & = &
\begin{array}{|c|c|}
\hline
100 & 60 \\ \hline
20 & 10 \\ \hline
\end{array}\ ,
\end{eqnarray*}
\end{minipage}
&
\begin{minipage}{100pt}
\begin{eqnarray*}
	T_2^{(2)} & = &
\begin{array}{|c|c|}
\hline
1 & 19 \\ \hline
5 & 20 \\ \hline
\end{array}\ ,
\end{eqnarray*}
\end{minipage}
&
\begin{minipage}{100pt}
\begin{eqnarray*}
	T_1^{(2)}+T_2^{(2)} & = &
\begin{array}{|c|c|}
\hline
101 & 79 \\ \hline
25 & 30 \\ \hline
\end{array}\ .
\end{eqnarray*}
\end{minipage} \\

\begin{minipage}{100pt}
\begin{eqnarray*}
	T_1^{(3)} & = &
\begin{array}{|c|c|}
\hline
100 & 60 \\ \hline
20 & 10 \\ \hline
\end{array}\ ,
\end{eqnarray*}
\end{minipage}
&
\begin{minipage}{100pt}
\begin{eqnarray*}
	T_2^{(3)} & = &
\begin{array}{|c|c|}
\hline
20 & 90 \\ \hline
5 & 20 \\ \hline
\end{array}\ ,
\end{eqnarray*}
\end{minipage}
&
\begin{minipage}{100pt}
\begin{eqnarray*}
	T_1^{(3)}+T_2^{(3)} & = &
\begin{array}{|c|c|}
\hline
120 & 150 \\ \hline
25 & 30 \\ \hline
\end{array}\ .
\end{eqnarray*}
\end{minipage} \\

\end{tabular}
\end{center}
Note that $(T_1^{(0)},T_2^{(0)})$ is an instance of positive SP, $(T_1^{(1)},T_2^{(1)})$ is an instance of positive alignment, 
$(T_1^{(2)},T_2^{(2)})$ is an instance of negative SP, and $(T_1^{(3)},T_2^{(3)})$ is an instance of negative aligment. 

For $k\geq 4$, define $(T_1^{(k)},T_2^{(k)})=(\delta T_1^{(k-4)},\delta T_2^{(k-4)})$, where $\delta=20$. Then 
$\mathcal{T}=(T_1^{(k)},T_2^{(k)})_{k\geq 0}$ is a monotonic sequence possessing the properties mentioned in the theorem.
\end{proof}

\end{document}